\documentclass[12pt]{article}
\title{ \Large Structural Properties of ${\cal R}_2$ Part I }
\author{\normalsize Timothy J. Carlson \\
 \normalsize The Ohio State University, Columbus, OH 43210 USA \\
\normalsize  email: carlson@math.ohio-state.edu \\
}

\date{\vspace{ -0.3in}}

\vspace{-.2in}
\usepackage{amsmath}    
\usepackage{amssymb}
\usepackage{graphicx}   
\usepackage{verbatim}   
\usepackage{color}      
\usepackage{subfigure}  
\usepackage{hyperref}   

\newcommand{\qed}{\hfill {$\square$}}
\newcommand{\blankline}{\vspace{4 mm}}
\newcommand{\halfblankline}{\vspace{2 mm}}

\newcommand{\bfA}{\mbox{${\bf A}$}}

\newcommand{\bfB}{\mbox{${\bf B}$}}

\newcommand{\calB}{\mbox{${\cal B}$}}

\newcommand{\calC}{\mbox{${\cal C}$}}

\newcommand{\calL}{\mbox{${\cal L}$}}

\newcommand{\calR}{\mbox{${\cal R}$}}

\newcommand{\calRrho}{\mbox{${\cal R}^\rho$}}
\newcommand{\calRtwo}{\mbox{${\cal R}_2$}}
\newcommand{\calRtwoalpha}{\mbox{${\cal R}_2^\alpha$}}
\newcommand{\calRtworho}{\mbox{${\cal R}_2^\rho$}}

\newcommand{\leqone}{\mbox{$\leq_1$}}
\newcommand{\leqonerho}{\mbox{$\leq_1^\rho$}}

\newcommand{\leqtwo}{\mbox{$\leq_2$}}
\newcommand{\leqtwodown}{\mbox{$\leq_2\downarrow$}}
\newcommand{\leqtworho}{\mbox{$\leq_2^\rho$}}
\newcommand{\leqtworhodown}{\mbox{$\leq_2^\rho\downarrow$}}
\newcommand{\leqtworhoup}{\mbox{$\leq_2^\rho\uparrow$}}
\newcommand{\leqtwoup}{\mbox{$\leq_2\uparrow$}}
\newcommand{\leqk}{\mbox{$\leq_k$}}

\newcommand{\lessone}{\mbox{$<_1$}}

\newcommand{\lesstwo}{\mbox{$<_2$}}

\newcommand{\precbfA}{\mbox{$\prec^{\bf A}$}}

\newcommand{\precbfB}{\mbox{$\prec^{\bf B}$}}

\newcommand{\preq}{\mbox{$\preceq$}}
\newcommand{\preqbfA}{\mbox{$\preceq^{\bf A}$}}

\newcommand{\preqone}{\mbox{$\preceq_1$}}

\newcommand{\preqtwo}{\mbox{$\preceq_2$}}

\newcommand{\unibfA}{\mbox{$|{\bf A}|$}}

\newcommand{\unibfB}{\mbox{$|{\bf B}|$}}

\newcommand{\unicalR}{\mbox{$|{\cal R}|$}}

\newcommand{\supbfA}{\mbox{$^{{\bf A}}$}}

\newtheorem{prop}{Proposition}[section]
\newtheorem{cor}[prop]{Corollary}
\newtheorem{thm}[prop]{Theorem}
\newtheorem{lem}[prop]{Lemma}

\newtheorem{dfn}[prop]{Definition}

\begin{document}

\maketitle


\noindent
{\small 
{\bf Abstract.} 
This is the first of two papers establishing structural properties of \calRtwo, the structure giving rise to pure patterns of resemblance of order two, which partially underly the results in [\ref{Ca09}] and [\ref{Ca16}] as well as other work in the area. 
}

\blankline

Letting $\leq$ denote the usual ordering of the class of ordinal numbers $ORD$ and using $\calB \preceq _{\Sigma _{n}} \calC$ to 
indicate that \calB\ is a $\Sigma _{n}$-elementary substructure of
\calC, define the structure  $\calRtwo \, = \, (ORD, \leq,\leqone,\leqtwo)$ by
inductively defining  the binary relations
$\leqone$ and $\leqtwo$ on $ORD$ so that 
\begin{center}
  $\alpha \leq_n \beta$ \ iff \ $(\alpha , \leq, \leqone,\leqtwo)
  \preceq _{\Sigma_{n}} (\beta ,  \leq,\leqone,\leqtwo)$
\end{center}
for $n=1,2$ and all ordinals $\alpha$ and $\beta$ (in other words, the restriction of \calRtwo\
to $\beta$ is defined by induction on $\beta$).

We will find it convenient to use an alternate definition for \calRtwo\ which was introduced in [\ref{Ca09}]. 
The equivalence of the two definitions will be established elsewhere.
Moreover, our investigation  of pure patterns of resemblance naturally lead to studying other variants of \calRtwo\ described below.

Some key properties of \calRtwo\ which are easily established are 
\begin{itemize}
 \item
 $\leq$ is a linear ordering of the ordinals.
\item
 \leqone\ and \leqtwo\ are partial orderings of the ordinals.
\item
 \leqtwo\ respects \leqone\ and \leqone\ respects $\leq$.
\end{itemize}
where for two partial orderings $\leq'$ and $\leq''$ on a set $X$ we say that $\leq''$ {\it respects} $\leq'$ if
$$x\leq'' y \ \ \ \Longrightarrow \ \ \ x\leq' y$$
for all $x,y\in X$ and
$$x\leq' y \leq' z \ {\rm and } \ x\leq'' z \ \ \ \Longrightarrow \ \ \ x\leq'' y$$
for all $x,y,z\in X$.

The main results of this paper, the various {\it recurrence theorems}, establish the recurrence of certain configurations as one proceeds through the ordinals.
Understanding this process for \calRtwo\ leads in a natural way to study slight variants of \calRtwo.
These results hinge on the cornerstone of the paper: the {\it Main Structural Lemma}.
While not difficult to prove, the Main Structural Lemma has far reaching consequences.

We define the {\it components} of $ORD$ with respect to \leqone\ to be the usual connectivity components of \leqone.
The fact that \leqone\ respects $\leq$ implies that the components  are intervals of the form $[\kappa,\lambda)$ where $\kappa$ and $\lambda$ are successive elements among the set of ordinals which are minimal in \leqone.
In fact, the components of \leqone\ are easily seen to be closed intervals.
These intervals are enumerated by $I_\alpha$ $(\alpha\in \theta_1)$.

The {\it First Recurrence Theorem for \leqone} says that the isomorphism type of $I_\alpha$ depends only on the last component of the additive normal form of $\alpha$ i.e. $I_{\beta+\omega^\alpha}\cong I_{\omega^\alpha}$.

The {\it Second Recurrence Theorem for \leqone} says 
$$(min(I_{\omega^\alpha}),max(I_{\omega^\alpha})]\cong (0,max(I_\alpha)]$$
This reduces the structure of $I_{\omega^\alpha}$ to that of earlier components in the case $\omega^\alpha$ is not an epsilon number.

For the remainder of the introduction, assume $\omega^\alpha$ is an epsilon number i.e. $\omega^\alpha=\alpha$.

Partition $I_\alpha$ into intervals as follows.
Let $X$ be the collection of $\nu\in I_\alpha$ such that $\nu$ is minimal with respect to \leqtwo\ and $\nu\leqone \tau$ whenever $\nu\leq \tau\in I_\alpha$. 
The {\it components} of  $I_\alpha$ with respect to \leqtwo\ are intervals of the form $[\nu,\mu)$ where $\nu$ and $\mu$ are two successive elements of $X$ along with $\{\tau \, | \, \nu\leq\tau\in I_\alpha\}$ when $\nu$ is the maximal element of $X$. 
Since \leqtwo\ respects \leqone, elements from different components are never related by \leqtwo. 
The components of $I_\alpha$ with respect to \leqtwo\ are enumerated by $J_\eta$ $(\eta\in\theta_2)$.

Let $\kappa$ be the least element of $I_\alpha$. 
The {\it Recurrence Theorem for Small Intervals} says that $(\mu,\mu+\kappa)\cong (0,\kappa)$ whenever $\mu$ is divisible by $\kappa$.
Moreover, $\nu\not\leq_2 \delta$ whenever $\nu\leq\mu<\delta<\mu+\kappa$.

The {\it First Recurrence Theorem for \leqtwo} says that the isomorphism type of $J_\eta$ depends only on the last component of the additive normal form of $\eta$ i.e $J_{\xi+\omega^\eta}\cong J_{\omega^\eta}$.

Section \ref{sp} contains background.

Section \ref{sv} contains the main definitions for the alternate versions of \calRtwo.

Section \ref{sbl} establishes general facts about the notions from Section \ref{sv}.

The Main Structural Lemma is presented in Section \ref{smsl}.

Section \ref{sfirstrone} establishes the First Recurrence Theorem for \leqone.

Section \ref{sis} develops the notion of {\it incompressible sets}, a weakened version of isominimal sets from [\ref{Ca09}] sufficient for the purposes of this paper.

Section \ref{ssecondrone} establishes the Second Recurrence Theorem for \leqone.

Section \ref{sfirstrtwo} establishes the Recurrence Theorem for Small Intervals and the First Recurrence Theorem for \leqtwo.


The sequel to this paper will contain the {\it Second Recurrence Theorem for \leqtwo}, an analogue of the Second Recurrence Theorem for 
\leqone\ which provides a characterization of the growth of the intervals $J_{\omega^\eta}$ as $\eta$ increases. 
This characterization hinges on the {\it Order Reduction Theorem} which describes $J_{\omega^\eta}$ in terms of an initial segment of a variant of \calRtwo\ denoted by \calRtwoalpha. 

We define \calRtwoalpha\ whenever $\alpha$ is an additively indecomposable ordinal.
The results of this paper for arbitrary $\alpha$ will only be necessary in the sequel to this paper. 
The present paper is self-contained if one ignores the variants and assume $\alpha=1$ throughout.
As will be clear, the case $\alpha=1$ is essentially \calRtwo.


\vspace{-4mm}

\section{ Preliminaries}
\label{sp}

As a general principle, we will omit parameters in defined notions when they are understood by the context. 

${\sf KP}$ will be used to denote Kripke-Platek set theory (see
[\ref{Ba75}] for background) and ${\sf KP}\omega$ is Kripke-Platek set theory
with the axiom of infinity. ${\sf KP}\omega$ is the base theory for the 
results in this paper. The theory ${\sf KP}\ell_{0}$ has an axiomatization consisting of the usual axiomatization
for ${\sf KP}\omega$ with $\Delta_{0}$-comprehension removed and an
additional axiom saying that every set is an element of an admissible set.
${\sf ZF}$ denotes Zermelo-Fraenkel set theory.

We will write $card(X)$ for the cardinality of a set $X$.
$ORD$ will denote
the class of ordinals with the usual ordering $\leq$ and arithmetic operations. 
$\omega$ is the least infinite ordinal and the 
elements of $\omega$ are natural numbers. 
 An ordinal is {\it additively indecomposable} if it is not $0$ and is closed under addition. 
The additively indecomposably ordinals can be characterized as the ordinals of the form $\omega^\beta$ for some $\beta$. 
For ordinals $\alpha$ and $\beta$, $\alpha$ {\it divides} $\beta$ if there exists and ordinal $\delta$ such that $\beta=\alpha\cdot\delta$.
We will also say either $\beta$ is divisible by $\alpha$ or $\beta$ is a multiple of $\alpha$ to mean $\alpha$ divides $\beta$.
$\beta$ is a {\it limit multiple} of $\alpha$ if $\beta=\alpha\cdot\delta$ for some limit ordinal $\delta$.
Similarly, $\beta$ is a {\it successor multiple} of $\alpha$ if $\beta=\alpha\cdot\delta$ for some successor ordinal $\delta$.

For $A$ a finite set of ordinals and $i$ less than the cardinality of $A$, $(A)_i$ is the $i^{th}$ element of $A$. 
We extend $<$ to finite sets of ordinals in the usual way e.g. $X<Y$ means that $\xi < \eta$ for all $\xi\in X$ and $\eta\in Y$. 
We will also write $X\not\leq_k Y$ to mean $\xi\not\leq_k\eta$ whenever $\xi\in X$ and $Y\in \eta$.

Contrary to standard practice, we will allow structures for a first-order language $\cal L$ to interpret the function symbols as partial operations on the universe which fail to give an interpretation to some constant symbols. 
In other words, we use the word ``structure'' to refer to what are called partial structures elsewhere (see [\ref{Be85}]).
We will write \unibfA\ for the universe of a structure \bfA. 
The definition of when a  term is defined in a structure is the natural one, proceeding from ``bottom up", as is the definition of  the value of the term in the structure. 
When $t$ is a term all of whose variables are among $v_1,\ldots,v_n$ and  $a_1,\ldots,a_n\in\unibfA$ we write $t(a_1,\ldots,a_n)\supbfA$ for the value of $t$ in \bfA, if it exists, when$v_1,\ldots,v_n$ are interpreted as $a_1,\ldots,a_n$ respectively. 
See the theory of partial terms in [\ref{Be85}] for details. 

When \bfA\ is a structure for $\cal L$ and $X\subseteq \unibfA$, the set {\it generated} in \bfA\ from $X$ is the smallest subset of \unibfA\ containing $X$ and closed under the interpretation of the function and constant symbols of $\cal L$ in \bfA.

We fix a special symbol $\preq$ which will be assumed to be
a 2-place relation symbol in every language in which it occurs.  Suppose \bfA\ is a structure for the first-order language $\cal L$ which includes
\preq.  If  the interpetation
of \preq\ in \bfA\ is a linear ordering of \unibfA\ we will say that
\bfA\ is a {\it linearly ordered structure}. 
If the interpretation of \preq\ in \bfA\ is a well ordering of \unibfA\ we will say that  \bfA\ is a {\it well ordered structure}.

Fix a linearly ordered structure \bfA.

If $X$ is 
a  nonempty subset of \unibfA, $max(X)$ will be the largest
element of $X$ and $min(X)$ will be the smallest element of $X$ if such elements
exist.

We will use standard interval notation e.g. for $a,b\in\unibfA$ let
$[a,b)\supbfA$ denote the set of all $x\in\unibfA$ such that
$a\preqbfA x \precbfA b$ (we write \precbfA\ for the strict part of the 
linear ordering \preqbfA). We will also use $\infty$ and $-\infty$ as interval endpoints in the usual ways e.g. write $(-\infty, a)\supbfA$ for the
set of all $x\in\unibfA$ such that $x\precbfA a$ and $[a,\infty)\supbfA$ for
the set of all $x\in\unibfA$ such that $a\preqbfA x$. 

Assume $X$ is  a subset of \unibfA. For $x$ an element of \unibfA\ which is not the least element of \bfA,
$x$ is a {\it limit point} of $X$ if $X$ intersects $(y,x)\supbfA$ whenever $y\precbfA x$ and $x$ is not minimal in \bfA. 
We say $x$ is a {\it limit element} of \bfA\ if it is a limit point of \unibfA.
$X$ is {\it topologically closed} in \bfA\ if every limit point of $X$ is an element of $X$.

If $y\precbfA x$ and $y$ is the largest element of $(-\infty,x)$ then $x$ is the {\it successor} of $y$.

Suppose $x$ is a limit element of \bfA. A subset $X$ of $(-\infty, x)$ is {\it cofinal} in $x$ if $x$ is a limit point of $X$. Similarly, a collection $\cal X$ of finite subsets of $(-\infty,x)$ is {\it cofinal} in $x$ if $\cal X$ contains a subset of $(y,x)\supbfA$ whenever $y\precbfA x$.

A subset $I$ of \unibfA\ is an {\it initial segment} of \bfA\ if
$x\in\ I$ whenever $x\preqbfA y$ for some $y\in I$. An initial segment $I$ of
\bfA\ is a {\it proper initial segment of \bfA} if $I\not = \unibfA$.
A substructure \bfB\ of \bfA\ is an {\it initial substructure of \bfA} if \unibfB\
is an initial segment of \bfA. \bfB\ is a {\it proper initial substructure of \bfA}
if the universe of \bfB\ is a proper initial substructure of \bfA.

An element $a$ of \unibfA\ is {\it decomposable} in \bfA\ if there is a function symbol $\sf f$  and $a_1,\ldots,a_n \precbfA a$ such that $a={\sf f}(a_1,\ldots,a_n)$ (we view constant symbols as 0-ary function symbols). If $a$ is not decomposable in \bfA\ then $a$ is {\it indecomposable} in \bfA. A substructure \bfB\ of \bfA\ is a {\it closed substructure} of \bfA\ if every element of \bfB\ which is indecomposable in \bfB\ is indecomposable in \bfA, or, equivalently, whenever $b\in\unibfB$ is decomposable in \bfA\ then there is a function symbol $\sf f$ and there are $b_1,\ldots,b_n\in\unibfB$ such that $b_1,\ldots,b_n\precbfB b$ and $b={\sf f}(b_1,\ldots,b_n)$. A subset of \unibfA\ is a {\it closed subset} of \bfA\ if it is the universe of a closed substructure of \bfA\ (note that we do not require a substructure of \bfA\ to be closed under the interpretations of the function symbols in \bfA\ so that any subset of \unibfA\ is the universe of a substructure of \bfA). 
Clearly, any set of indecomposable elements is closed and the closed subsets are closed under initial segments and arbitrary unions.
Also, if \bfA\ is a well ordered structure for \calL\ then every finite subset of \unibfA\ is contained in a finite closed subset of \unibfA. 

The closed sets in variants of \calRtwo\ will be particularly simple. In particular, there are no functions in \calRtwo\ itself implying every ordinal is indecomposable and every set is closed.

\section{Variants of \calRtwo}
\label{sv}

We will assume $\rho$ is an arbitrary additively indecomposable ordinal for the rest of the paper.

We will introduce structures similar to  \calRtwo\ which are useful in our analysis.

\begin{dfn}
For $0\leq\xi<\rho$ define $f^\rho_\xi:ORD\rightarrow ORD$ so that 
\begin{equation}
f^\rho_\xi(\rho\cdot \zeta + \epsilon) =
\begin{cases}
\rho\cdot\zeta + \xi & \text{if $\epsilon = 0$} \\
\rho\cdot\zeta & \text{if $\epsilon \not=0$}
\end{cases}
\end{equation}
for all ordinals $\zeta$ and $\epsilon$ with $\epsilon< \rho$.
${\cal R}^\rho$ is the structure 
 $$(ORD,f_\xi^\rho \ (0\leq\xi < \rho),\leq)$$
 \end{dfn}
 
 \halfblankline
 
\calRrho\ is an EM structure as defined in [\ref{Ca09}]. 

Notice that $f^\rho_0(\rho\cdot\zeta +\xi)=\rho\cdot\zeta$ whenever $\xi<\rho$.
Therefore, $f^\rho_0(\alpha)=\alpha$ iff $\alpha$ is divisible by $\rho$.

 \begin{lem}
 \label{lvA} \
 \begin{enumerate}
 \item
 The indecomposable ordinals of \calRrho\ are the ordinals which are divisible by $\rho$ i.e. ordinals of the form $\rho\cdot\zeta$.
 \item
 An ordinal $\alpha$ is indecomposable in \calRrho\ iff $f^\rho_0(\alpha)=\alpha$.
 \item
 A set of ordinals $X$ is closed in \calRrho\ iff $\sigma\in X$ whenever $\sigma$ is divisible by $\rho$ and $\sigma+\xi\in X$  for some $\xi<\rho$ i.e. $X$ is closed under $f_0^\rho$.
 \item
 Any union  or intersection of closed sets is closed.
  \item
An initial segment of a closed set is closed.
 \item
 For all ordinals $\alpha$, $\alpha$ is closed.
 \item
 For all ordinals $\alpha$, $[\alpha,\infty)$ is closed iff $\alpha$ is divisible by $\rho$.
 \item 
 If $X$ is closed and $\alpha$ is divisible by $\rho$ then $X\cap[\alpha,\infty)$ is closed.
 \end{enumerate}
 \end{lem}
 {\bf Proof.} 
Straigtforward.

For example, to prove part 3, suppose $X$ is a set of ordinals.

First, assume $X$ is closed in \calRrho.
Also, suppose $\sigma+\xi\in X$ where $\sigma$ is divisible by $\rho$ and $\xi<\rho$.
We must show $\sigma\in X$.
This is trivial if $\xi=0$, so we may assume $0<\xi$. 
Since $\sigma+\xi=f_\xi(\sigma)$, $\sigma+\xi$ is decomposable. 
Therefore, there exists $\eta<\rho$ and $\tau \in X$ such that $\tau<\sigma+\xi$ such that $f_\eta(\tau)=\sigma+\xi$.
Clearly, $\tau=\sigma$ and $\eta=\xi$. 

Now, assume $\sigma\in X$ whenever $\sigma$ is divisible by $\rho$ and $\sigma+\xi\in X$ for some $\xi<\rho$.
To show $X$ is closed, assume $\alpha\in X$ is decomposable.
There exists $\xi<\rho$ and $\sigma<\alpha$ such that $f_\xi(\sigma)=\alpha$.
Since $f_\xi(\sigma)>\sigma$, $\sigma$ is divisible by $\rho$ and $f_\xi(\sigma)=\sigma+\xi$.
By assumption, $\sigma\in X$.
 \qed
 
\begin{dfn}
For $\alpha$ an ordinal, $rem^\rho(\alpha)$, the {\bf remainder} of $\alpha$ with respect to $\rho$, is the unique $\varepsilon<\rho$ such that $\alpha=\rho\cdot\delta + \varepsilon$ for some $\delta$.
\end{dfn}
 
 As usual, when $\rho$ is clear from the context, we will simply write $rem$ for $rem^\rho$.


\begin{lem}
\label{lvB}
Assume $X$ is a closed set of ordinals and $h$ is an order preserving function which maps $X$ into $ORD$.   
The following are equivalent.
\begin{enumerate}
\item
$h$ is an embedding of $X$, as a substructure of \calRrho, into \calRrho.
\item
For all ordinals $\sigma$ and $\xi$, if $\sigma$ is divisible by $\rho$, $\xi<\rho$ and $\sigma+\rho\in X$ then $h(\sigma+\xi)=h(\sigma)+\xi$ and $h(\sigma)$ is divisible by $\rho$.
\item
 The range of $h$ is closed and $rem(\alpha)=rem(h(\alpha))$ for all $\alpha\in X$.
 \end{enumerate}
\end{lem}
{\bf Proof.}
$(1\Rightarrow 2)$ Follows from the definition of $f_\xi$ and parts 1 and 2 of the previous lemma.

$(2\Rightarrow 3)$ Assume part 2.

We will use part 3 of the previous lemma to show $h[X]$ is closed.
Suppose $\sigma+\xi\in h[X]$ where $\sigma$ is divisible by $\rho$ and $\xi<\rho$.
There exist $\alpha\in X$ such that $\sigma+\xi=h(\alpha)$.
There exist $\tau$ and $\eta$ such that $\tau$ is divisible by $\rho$, $\eta<\rho$ and $\alpha=\tau+\eta$.
Since $X$ is closed, $\tau\in X$.
By part 2, $\sigma+\xi=h(\tau)+\eta$ and $h(\tau)$ is divisible by $\rho$.
Clearly, $\sigma=h(\tau)$.

Now suppose $\alpha\in X$. 
We will show $rem(h(\alpha))=rem(\alpha)$.
Let $\xi=rem(\alpha)$.
There exists $\sigma$ such that $\sigma$ is divisible by $\rho$ and $\alpha=\sigma+\xi$.
By part 2, $h(\alpha)=h(\sigma)+\xi$ and $h(\sigma)$ is divisible by $\rho$.
This implies $rem(h(\alpha))=\xi$.

$(3 \Rightarrow 1)$
Notice that for any ordinals $\alpha$, $\beta$ and $\xi$ with $\xi<\rho$, $f_\xi(\alpha)=\beta$ iff either 
\begin{itemize}
\item[]
$rem(\beta)=\xi$ and $\alpha$ is the largest  $\gamma$ such that $\gamma\leq\beta$ and $rem(\gamma)=0$
\end{itemize}
or
\begin{itemize}
\item[]
$rem(\alpha)\not=0$ and $\beta$ is the largest  $\gamma$ such that $\gamma\leq\alpha$ and $rem(\gamma)=0$.
\end{itemize}
The implication $(3\Rightarrow 1)$ follows from the observation that if $Z$ is a closed set of ordinals and $\delta,\varepsilon\in Z$ then $\delta$ is the largest $\gamma$ such that $\gamma\leq\varepsilon$ and $rem(\gamma)=0$ iff $\delta$ is the largest $\gamma\in Z$ such that $\gamma\leq\varepsilon$ and $rem(\gamma)=0$. 
\qed

\begin{lem}
\label{lvBB}
Assume $X$ is a closed set of ordinals.
If $h$ is an order preserving map of the indecomposable ordinals in $X$ into the class of indecomposable ordinals there is a unique extension $h^+$ to $X$ which is an embedding of $X$, as a substructure of \calRrho, into \calRrho.
\end{lem}
{\bf Proof.}
By the previous lemma.
\qed

\begin{lem}
\label{lvB1}
\begin{enumerate}
\item
For all $\alpha\in ORD$, $[\rho\cdot\alpha,\infty)\cong ORD$.
\item
For all $\alpha,\xi\in ORD$, if $0<\xi<\rho$ then $[\alpha+\xi,\infty)\cong[1,\infty)$.
\end{enumerate}
\end{lem}
{\bf Proof.}
Part 1 follows from the Lemma \ref{lvB}. 

For part 2, notice part 1 implies that $[\alpha+\rho,\infty)\cong[\rho,\infty)$ since $\rho$ is additively indecomposable. Now notice $[\alpha+\xi,\alpha+\rho)$ and $[1,\rho)$ both have length $\rho$ and nothing in either set is either in the range or domain of the restriction any $f_\xi$.
\qed

\begin{lem}
\label{lvC}
Assume $X$ is a closed set of ordinals and $h$ is an embedding of $X$, as a substructure of \calRrho, into \calRrho. 
If $\alpha$ is the least ordinal moved by $h$ then $\alpha$ is divisible by $\rho$.
\end{lem}
{\bf Proof.}
Assume $\alpha$ is the least element moved by $h$.
There are $\zeta$ and $\xi<\rho$ such that $\alpha=\rho\cdot\zeta+\xi$.
Argue by contradiction and assume $\xi\not=0$.
Since $X$ is closed, $\rho\cdot\zeta\in X$.
Since $\rho\cdot\zeta<\alpha$, $h(\rho\cdot\alpha)=\rho\cdot\alpha$.
Since $f^\rho_\xi(\rho\cdot\zeta)=\alpha$, 
$$h(\alpha)=h(f^\rho_\xi(\rho\cdot\zeta))=f^\rho_\xi(h(\rho\cdot\zeta))=f^\rho_\xi(\rho\cdot\zeta)=\alpha$$
-- contradiction.
\qed

\begin{dfn}
Assume \calR\ is a structure of the form 
$$(ORD, f^\rho_\xi \ (\xi<\rho), \leq, \preqone, \preqtwo)$$
or the restriction of such a structure to some ordinal.
 Also assume $X$ is a closed set of ordinals (with respect to \calRrho) which is a subset of \unicalR.
A function $h:X\rightarrow \unicalR$ is a {\bf covering}
of $X$ in \calR\ if 
\begin{enumerate}
 \item $h$ is an embedding of  $X$ as a substructure of \calRrho\ into
\calRrho.
\item  For $\alpha,\beta\in X$ and $k=1,2$. 
 $$\alpha\preq^\rho_k\beta \ \ \  \Longrightarrow \ \ \ h(\alpha)\preq^\rho_k h(\beta)$$
 \end{enumerate}
 A set $Y$ of ordinals is a {\bf covering} of $X$ in \calR\ if there is a function $h$ which is a covering of $X$ in \calR\ with range $Y$.
\end{dfn}

We will say that a set of ordinals $Y$ is a covering of $X$  if there is a covering of $X$ in \calR\ with range $Y$ and \calR\ is clear from the context. 

\begin{lem} 
Assume \calR\ is as above. 
\begin{enumerate}
\item
The composition of coverings is a covering. 
\item
If $X$ is a covering of $Y$ and $Y$ is a covering of $Z$ then $X$ is a covering of $Z$.
\item The range of a covering is closed.
\end{enumerate}
\end{lem}
{\bf Proof.} Clear.
\qed

 \begin{dfn}
We define the structure 
$$\calRtworho=(ORD, f^\rho_\xi (0\leq \xi < \rho),\leq,\leqonerho, \leqtworho)$$
 so that for $\beta\in ORD$
 the following recursive clauses hold:
 \begin{itemize}
 \item   For all $\alpha\in ORD$
 $$\alpha \leqonerho \beta$$ 
 iff 
 \begin{itemize}
 \item[]
$\alpha\leq\beta$
\end{itemize}
and
\begin{itemize}
\item[]
 for any finite  $X\subseteq \alpha$ and finite $Y\subseteq [\alpha,\beta)$ with $X\cup Y$ closed, there is a finite $\tilde{Y}\subseteq \alpha$ such that
  \begin{itemize}
   \item $X < \tilde{Y}$
   \item $X\cup \tilde{Y}$ is a covering of  $X \cup Y$
  \end{itemize}
  \end{itemize}
 \item   For all $\alpha\in ORD$
 $$\alpha \leqtworho \beta$$ 
 iff  
 \begin{itemize}
 \item[] $\alpha\leq \beta$
 \end{itemize}
 and
 \begin{enumerate}
 \item[] for any finite sets $X$ and  $Y$ below $\alpha$ with $X$ and $Y$  closed and $X<Y$,
  if there are cofinally many sets $\tilde{Y}$ below $\alpha$
   such that  $X\cup \tilde{Y}$ is a covering
   of $ X\cup Y$
  then there are cofinally many 
   sets $\tilde{Y}$ below $\beta$ such that $X\cup \tilde{Y}$ is a covering of $X\cup Y$.
\end{enumerate}
and
\begin{itemize}
 \item[] for any finite $X\subseteq \alpha$ and finite $Y\subseteq [\alpha,\beta)$ with $X\cup Y$ closed, there is a finite $\tilde{Y}\subseteq \alpha$ such that
  \begin{itemize}
   \item $X < \tilde{Y}$
   \item $X\cup \tilde{Y}$ is a covering of  $X \cup Y$
   \item For any $i< card(Y)$, if $(Y)_i\leqonerho \beta$   then $(\tilde{Y})_i\leqonerho \alpha$
    (recall $(Y)_i$ is the $i^{\rm th}$ element of $Y$).
  \end{itemize}
 \end{itemize}
\end{itemize}

\end{dfn}

\halfblankline

Notice that
\begin{itemize}
 \item Whether $\alpha\leqonerho\beta$ holds depends only on the restriction of \calRtworho\ to $\beta$.
 \item Whether $\alpha\leqtworho\beta$ holds depends only on the restriction of \calRtworho\ to $\beta$ and the collection of $\alpha$ such that $\alpha\leqonerho\beta$.
 \end{itemize}

We will find it convenient to have separate notation for the properties given by the final two clauses in the definition of \leqtworho.

\begin{dfn}
Assume $\alpha,\beta\in ORD$ and $\alpha\leq\beta$. Define
$$\alpha\leqtworhoup \beta$$
iff
 \begin{enumerate}
 \item[] for any finite sets $X$ and  $Y$ below $\alpha$ with $X\cup Y$ closed and $X<Y$,
  if there are cofinally many sets $\tilde{Y}$ below $\alpha$
   such that  $X\cup \tilde{Y}$ is a covering
   of $ X\cup Y$
  then there are cofinally many 
   sets $\tilde{Y}$ below $\beta$ such that $X\cup \tilde{Y}$ is a covering of $X\cup Y$.
\end{enumerate}
Also, define
$$\alpha\leqtworhodown\beta$$
iff
\begin{itemize}
 \item[] for any finite $X\subseteq \alpha$ and finite $Y\subseteq [\alpha,\beta)$ with $X\cup Y$ closed, there is a finite $\tilde{Y}\subseteq \alpha$ such that
  \begin{itemize}
   \item $X < \tilde{Y}$
   \item $X\cup \tilde{Y}$ is a covering of  $X \cup Y$
   \item For any $i< card(Y)$, if $(Y)_i\leqonerho \beta$   then $(\tilde{Y})_i\leqonerho \alpha$
    (recall $(Y)_i$ is the $i^{\rm th}$ element of $Y$).
  \end{itemize}
 \end{itemize}
\end{dfn}

The definition of \leqtworhoup\ has the peculiarity that $\alpha\leqtworhoup\beta$ vacuously whenever $\alpha$ is not a limit multiple of $\rho$ i.e. not of the form $\rho\cdot\lambda$ for some limit $\lambda$.
On the other hand, we will only be interested in the case when $\alpha$ is a limit multiple of $\rho$.

Our definition of \calRtworho\ is a special case of the definition of \calRtwo\ in  Definition 5.4 of [\ref{Ca09}].

Notice that $f^1_0$ is the identity function. The structure obtained by removing $f^1_0$ from ${\cal R}^1_2$ is the same as \calRtwo\  as  defined in the introduction  using relations of partial elementarity (this will be established elsewhere). Hence, the results for general \calRtworho\ in the remainder of the paper include the version of \calRtwo\ from the introduction as a special case.

\section{Basic Lemmas}
\label{sbl}

Recall that  $\rho$ is an arbitrary additively indecomposable ordinal.

Since $\rho$ is fixed we will often omit $\rho$ as a parameter in notation e.g. we will write $f_\xi$, \leqone, \leqtwodown, \leqtwoup\ and \leqtwo\ for $f^\rho_\xi$, \leqonerho, \leqtworhodown, \leqtworhoup\ and  \leqtworho\ respectively.

Following the usual convention, we will write $\alpha<^\rho_k\beta$ when $\alpha\leq^\rho_k \beta$ and $\alpha\not=\beta$.


\begin{lem}
\label{lblD1}
If $\alpha\lessone \beta $ then $\alpha$ is a limit multiple of $\rho$.
\end{lem}
{\bf Proof.}
There are $\lambda$ and $\eta$ with $\eta<\rho$ such that $\alpha=\rho\cdot\lambda + \eta$. 

We claim $\alpha=\rho\cdot\lambda$.
Argue by contradiction and assume $0<\eta$.
This implies $\rho\cdot\lambda<\alpha$.
Let $X=\{\rho\cdot\lambda\}$ and $Y=\{\alpha\}$. 
Notice that $X\cup Y$ is closed.
Since $\alpha<_1\beta$, there is $\tilde{\alpha}<\alpha$ with $\rho\cdot\lambda<\tilde{\alpha}$ such that $\{\rho\cdot\lambda,\tilde{\alpha}\}$ is a covering of $\{\rho\cdot\lambda,\alpha\}$.
Since $f_\eta(\rho\cdot\lambda)=\alpha$, this implies  $f_\eta(\rho\cdot\lambda)=\tilde{\alpha}$ -- contradiction.

By letting $X=\emptyset$ and $Y=\{\alpha\}$ in the definition of \leqone\ we see there exists $\tilde{\alpha}<\alpha$. 
Therefore, $\alpha\not=0$.

To see that $\lambda$ is a limit ordinal, it suffices to show there are cofinally many ordinals below $\alpha$ which are divisible by $\rho$.

Suppose $\alpha'<\alpha$. 
There are $\lambda'$ and $\eta'$ with $\eta'<\rho$ such that $\alpha'=\rho\cdot\lambda' + \eta'$.
Let $X=\{\rho\cdot\lambda',\alpha'\}$ and $Y=\{\alpha\}$.
By part 2 of Lemma \ref{lvA}, $X$ and $Y$ are closed.
Since $\alpha<_1\beta$, there is $\tilde{\alpha}<\alpha$ with $\alpha'<\tilde{\alpha}$ such that  $\{\rho\cdot\lambda',\alpha',\tilde{\alpha}\}$ is a covering of $\{\rho\cdot\lambda',\alpha',\alpha\}$.
Since $f_0(\alpha)=\alpha$, $f_0(\tilde{\alpha})=\tilde{\alpha}$.
Therefore, $\tilde{\alpha}$ is divisible by $\rho$.
\qed

\begin{lem}
Assume   $\alpha<\beta$. 
\begin{enumerate}
\item
$\alpha\leqone\beta$ iff $\alpha$ is divisible by $\rho$ and the definition of $\alpha\leqonerho\beta$ holds with``$X\cup Y$ closed" replaced by ``$X$ and $Y$ closed".
\item
$\alpha\leqtwodown\beta$ iff $\alpha$ is divisible by $\rho$ and the definition of $\alpha\leqtworhodown\beta$ holds with ``$X\cup Y$ closed" replaced by ``$X$ and $Y$ closed".
\end{enumerate}
\end{lem}
{\bf Proof.}
By the previous lemma and parts 4, 5 and 8 of Lemma \ref{lvA}.
\qed

\blankline

We will use the previous lemma implicitly when checking whether $\alpha\leqone\beta$ or $\alpha\leqtwodown\beta$ for the rest of the paper.

The next two lemmas provide a useful alternative characterization of $\leqtworhoup$.

\begin{lem}
\label{lblD2}
Assume $n\in\omega$ and $R$ is a finite subset of $\rho$.
If $\cal Y$ is the collection of closed sets of ordinals $Y$ such that
\begin{itemize}
\item
the cardinality of $Y$ is at most $n$ and 
\item
 $rem[Y]\subseteq R$ 
\end{itemize}
then there are only finitely many elements of $\cal Y$ up to isomorphism as substructures of \calRtworho.
\end{lem}
{\bf Proof.} By Lemma \ref{lvB}, two finite closed sets are isomorphic as subsets of \calRrho\ iff they have the same cardinality and the order preserving map between them preserves remainders with respect to $\rho$. 
Hence, there are only finitely many elements of $\cal Y$ up to isomorphism as substructures of \calRrho.
Clearly, among any subcollection of $\cal Y$ whose elements are pairwise isomorphic as substructures of \calRrho, there are only finitely many isomorphism types as substructures of \calRtworho.
\qed

\begin{lem}
\label{lblD3}
Assume $\alpha,\beta\in ORD$. 
$$\alpha\leqtworhoup \beta$$
 iff 
 \begin{itemize}
\item[]
for any finite closed $X\subseteq \alpha$, for any family $\cal Y$ of nonempty finite closed sets which are cofinal in $\alpha$ with $X<Y$ for all $Y\in \cal Y$ and  $X\cup Y_1\cong X\cup Y_2$ whenever $Y_1,Y_2\in \cal Y$, and for any $\beta'<\beta$, there exists $Y\in \cal Y$  and $\tilde{Y}$ with $\beta'<\tilde{Y}<\beta$ such that $X\cup\tilde{Y}$ is a covering of $X\cup Y$.
\end{itemize}
\end{lem}
{\bf Proof.}
Straightforward using an easier version of the previous lemma in which one assumes the elements of $\cal Y$ are pairwise isomorphic as substructures of \calRrho.
\qed

\blankline

We could have just as well replaced ``there exists $Y\in \cal Y$ and $\tilde{Y}$" by ``for all $Y\in \cal Y$ there exists $\tilde{Y}$" in the right hand side of the equivalence.
The form above is slightly more convenient to apply.

\begin{lem}
\label{lblA}
\ 
\begin{enumerate}
 \item If $\alpha \leqone \beta$ and   $X\subseteq \alpha$ and $Y \subseteq [\alpha,\beta)$ are finite with $X$ and $Y$ closed  then there are cofinally many $\tilde{Y}$ below $\alpha$ such that $X\cup \tilde{Y}$ is a covering of $X\cup Y$.
  \item If $\alpha \leqtwodown \beta$ and  $X\subseteq \alpha$ and $Y \subseteq [\alpha,\beta)$ are finite with $X$ and $Y$ closed then there are cofinally many $\tilde{Y}$ below $\alpha$ such that $X\cup \tilde{Y}$ is a covering of $X\cup Y$ and $(\tilde{Y})_i \leqone \alpha$ whenever $i<card(Y)$ and $(Y)_i\leqone \beta$.
  \item Assume $\alpha,\beta_1,\beta_2\in ORD$ with $\alpha\leq \beta_1, \beta_2$.
   If there are $\gamma_i \in [\alpha,\beta_i)$ for $i=1,2$ such that $\alpha \cup (\gamma_1,\beta_1) \cong \alpha \cup (\gamma_2,\beta_2)$ then $\alpha \leqtwoup \beta_1$ iff $\alpha \leqtwoup \beta_2$.
 \end{enumerate}
 \end{lem}
 {\bf Proof.}
For part 1, assume $\alpha'< \alpha$. We need to show there exists $\tilde{Y}$ such that $\alpha'<\tilde{Y}<\alpha$ and $X\cup\tilde{Y}$ is a covering of $X\cup Y$. 
By  Lemma \ref{lblD1}, we may assume $\alpha'$ is divisible by $\rho$ and $X<\alpha'$.

Since $\alpha \leqone \beta$, there exists $\tilde{Y}\subseteq \alpha$ such that $X\cup\{\alpha'\}< \tilde{Y}$ and $X\cup\{\alpha'\} \cup \tilde{Y}$ is a covering of $X\cup\{\alpha'\}\cup Y$. $X\cup\tilde{Y}$ is a covering of $X\cup Y$.

The proof of part 2 is similar to the proof of part 1.

Part 3 is immediate.
 \qed

 \begin{lem}
 \label{lblD}
  \begin{enumerate}
   \item If $\alpha\lesstwo \beta$ then $\beta$ is a limit multiple of $\rho$ and there are cofinally many $\gamma$ below $\alpha$ such that $\gamma\lessone\alpha$.
   \item If $\alpha\lesstwo\beta\lesstwo\gamma$ then the collection of $ \delta$ such that $\alpha\lesstwo\delta$ is cofinal in $\gamma$.
   \end{enumerate}
\end{lem}
{\bf Proof.}
For part 1, notice that if $\gamma$ is divisible by $\rho$ (equivalently, $f_0(\gamma)=\gamma$) then $\{\gamma\}\cong\{\delta\}$ iff $\delta$ is divisible by $\rho$.
Since $\alpha\leqtwoup\beta$, Lemma \ref{lblD1} implies the collection of ordinals divisible by $\rho$ is cofinal in $\beta$.
Therefore, $\beta$ is a limit multiple of $\rho$.

By letting $X=\emptyset$ and $Y=\{\alpha\}$ in the definition of \leqtwodown, we see there are cofinally many $\gamma$ below $\alpha$ with $\gamma<_1\alpha$ by part 2 of Lemma \ref{lblA}.

For part 2, assume $\alpha\lesstwo\beta\lesstwo\gamma$.
By Lemma \ref{lblD1}, $\alpha$ and $\beta$ are divisible by $\rho$.
Hence, $\{\alpha\}$ and $\{\beta\}$ are closed.
Since $\beta\leqone\gamma$, there are cofinally many $\beta'$ below $\beta$ such that $\{\alpha,\beta'\}$ is a covering of $\{\alpha,\beta\}$. 
By Lemma \ref{lblD2}, there is a cofinal set of these $\beta'$ for which the isomorphism type of $\{\alpha,\beta'\}$ is fixed.
Choose $\beta_0$ to be one element of this set.
Since $\beta<_2\gamma$, there are cofinally many $\beta'$ below $\gamma$ such that $\{\alpha,\beta'\}$ is a covering of $\{\alpha,\beta_0\}$ and, hence, a covering of $\{\alpha,\beta\}$.
Any such $\beta'$ has the property that $\alpha\leqtwo\beta'$.
\qed

\begin{lem} 
\label{lblB}
\begin{enumerate}
 \item \leqone\ respects $\leq$.
  \item  \leqtwo\ respects \leqone.
 \item  \leqone\ is a partial ordering of $ORD$.
 \item \leqtwodown, \leqtwoup\ and \leqtwo\ are partial orderings of $ORD$.
\end{enumerate}
\end{lem}
{\bf Proof.} 
Part 1 is clear.

For part 2,  assume $\alpha\leqone \beta\leqone \gamma$ and $\alpha\leqtwo \gamma$. 
We will show $\alpha\leqtwo \beta$.
We may assume $\alpha< \beta< \gamma$ since the proof is trivial otherwise.

To show $\alpha\leqtwodown \beta$, assume $X\subseteq \alpha$ and $Y\subseteq [\alpha,\beta)$ are finite with $X$ and $Y$ closed. Since $\alpha \leqtwodown \gamma$, there exists $\tilde{Y}\subseteq \alpha$ such that $X< \tilde{Y}$, $X\cup \tilde{Y}$ is a covering of $X\cup Y$ and $(\tilde{Y})_i\leqone \alpha$ whenever $(Y)_i\leqone \gamma$. Since $\beta\leqone \gamma$, $(\tilde{Y})_i\leqone \alpha$ whenever $(\tilde{Y})_i\leqone \beta$.

To show $\alpha\leqtwoup \beta$, assume $X$ and $Y$  are finite subsets of $\alpha$ with $X$ and $Y$ closed and $X<Y$ such that there are cofinally many $\tilde{Y}$ below $\alpha$ such that $X\cup \tilde{Y}$ is a covering of $X\cup Y$. Since $\alpha\leqtwoup \gamma$, there exists $Y^*\subseteq [\beta,\gamma)$ such that $X\cup Y^*$ is a covering of $X\cup Y$. Since $\beta\leqone \gamma$, part 1 of Lemma \ref{lblA} implies that there are cofinally many $\tilde{Y}$ below $\beta$ such that $X\cup \tilde{Y}$ is a covering of $X\cup Y^*$. 

The proof of part 3 is similar and easier than that of part 4, so we will omit it.

The relations $\leqtwodown$, $\leqtwoup$ and $\leqtwo$ are clearly reflexive.
Since $\leqtwodown$, $\leqtwoup$ and $\leqtwo$ are contained in $\leq$, they are antisymmetric. The transitivity of $\leqtwo$ follows from that of \leqtwodown\ and \leqtwoup.

To see that \leqtwodown\ is transitive, assume $\alpha \leqtwodown \beta \leqtwodown \gamma$. 
To show $\alpha \leqtwodown \gamma$, assume $X\subseteq \alpha$ and $Y\subseteq [\alpha,\gamma)$ are finite with $X$ and $Y$ closed.  
Without loss of generality, assume $\alpha\in Y$.
Since $\beta\leqtwodown\gamma$, a simple argument shows there exists $Y^*\subseteq [\alpha,\beta)$ such that $X\cup Y^*$ is a covering of $X\cup Y$ and $(Y^*)_i\leqone \beta$ whenever $(Y)_i\leqone \gamma$. 
The desired $\tilde{Y}$ can be obtained using the assumption that $\alpha\leqtwodown \beta$.

The transitivity of \leqtwoup\ is straightforward.
\qed

\blankline

The lemma implies that $\leq_k$ is a forest on $ORD$ for $k=1,2$ \  i.e. a partial ordering in which the precedessors of any element are linearly ordered.

\begin{lem}
 \
 \label{lblE}
\begin{enumerate}
 \item Assume $0< \alpha <  \beta$ and $k\in\{1,2\}$. If for all $\alpha' <  \alpha$ and $\beta'< \beta$ there exist $\alpha''\in (\alpha',\alpha]$ and $\beta''\in (\beta',\beta]$ such that $\alpha'' \leq_k \beta''$ then $\alpha \leq_k \beta$.
 \item For $\alpha\in ORD$, the collection of $\beta$ in $ORD$ with $\alpha\leqone \beta$ is a topologically closed interval.
\item For $\alpha\in ORD$, the collection of $\beta\in ORD$ with $\alpha\leqtwo \beta$ is topologically closed.
\item For $\beta\in ORD$ and $k=1,2$, the collection of $\alpha\in ORD$ with $\alpha\leq_k \beta$ is topologically closed.
\end{enumerate}
\end{lem} 
{\bf Proof.}
For part 1, assume for all $\alpha' <  \alpha$ and $\beta'< \beta$ there exist $\alpha''\in (\alpha',\alpha]$ and $\beta''\in (\beta',\beta]$ such that $\alpha'' \leq_k \beta''$.

First assume $k=1$. To show $\alpha\leqone \beta$, suppose $X\subseteq \alpha$ and $Y\subseteq [\alpha,\beta)$ are finite with $X\cup Y$ closed. There are $\alpha''\leq \alpha$ and $\beta''\leq \beta$ such that $X< \alpha''$, $Y< \beta''$ and $\alpha''\leqone \beta''$. Since $\alpha''\leqone \beta''$, there exists $\tilde{Y}\subseteq \alpha''$ such that $X< \tilde{Y}$ and $X\cup\tilde{Y}$ is a covering of $X\cup Y$.

Now assume $k=2$ in part 1. 
The proof that $\alpha\leqtwodown \beta$ is similar to the proof for $k=1$. 
To show that $\alpha\leqtwoup \beta$, assume $X$ and $Y$ are finite   subsets of $\alpha$  with $X\cup Y$ closed such that there are cofinally many $\tilde{Y}$ below $\alpha$ such that $X\cup \tilde{Y}$ is a covering of $X\cup Y$. 
To show there cofinally many such $\tilde{Y}$ below $\beta$, assume $\beta'< \beta$. 
There are $\alpha''\leq \alpha$ and $\beta''\leq \beta$ such that $X< \alpha''$, $\beta'< \beta''$ and $\alpha''\leqtwo \beta''$. 
Since $\alpha''\leqone\beta''$, $\alpha''\leqone\alpha$.
Notice that there are cofinally many $\tilde{Y}$ below $\alpha''$ such that $X\cup \tilde{Y}$ is a covering of $X\cup Y$. 
This is an assumption if $\alpha''=\alpha$ and follows from part 1 of Lemma \ref{lblA} otherwise since $\alpha''\leqone\alpha$. Since $\alpha''\leqtwoup \beta''$, there are cofinally many $\tilde{Y}$ below $\beta''$ such that $X\cup \tilde{Y}$ is a covering of $X\cup Y$. In particular, there exists $\tilde{Y}$ with $\beta'<  \tilde{Y}< \beta''$ such that $X\cup \tilde{Y}$ is a covering of $X\cup Y$.

Parts 2-4 follow immediately from part 1.
\qed

\blankline

 The following lemma will be useful in establishing relations of the form $\alpha\leqone\beta$ and $\alpha\leqtwodown\beta$.
 
\begin{lem}
\label{lblC}
 Assume $X< Y_1< Y_2$ are finite closed sets of ordinals and let $X'$ be the set of $\xi\in X$ such that there is $\eta \in Y_2$ with $\xi\leqtwo\eta$. If $X'\cup Y_1$ is a covering of $X'\cup Y_2$ then $X\cup Y_1$ is a covering of $X\cup Y_2$.
\end{lem}
{\bf Proof.}
The lemma is trivial if the $Y_i$ are empty, so we may assume both $Y_1$ and $Y_2$ are nonempty.

Assume $h$ is a covering of $X'\cup Y_2$ onto $X'\cup Y_1$.
Since $Y_1<Y_2$, the least element of $Y_2$ must be divisible by $\rho$ by Lemma \ref{lvC}.
Extend $h$ to a function $h^+$ with domain $X\cup Y_2$ such that $h^+(\alpha)=\alpha$ for all $\alpha\in X$.
$h^+$ is an embedding by Lemma \ref{lvB}.

To see that $h^+$ is a covering, assume $\alpha,\beta\in X$ and $\alpha\leqk \beta$.
If $\alpha,\beta\in X$, $\alpha,\beta\in X'\cup Y_2$ or $\alpha=\beta$ then $hd^+(\alpha)\leqk h^+(\beta)$ is clear.
Consider the remaining case where $\alpha\in X-X'$ and $\beta\in Y_2$.
Since $\alpha\not\in X'$, $k=1$ and $\alpha\leqone\beta$.
Since $h(\alpha)=\alpha$ and $h^+(\beta)=h(\beta)$, it suffices to show $\alpha\leqone h(\beta)$.
Since $h(\beta)\in Y_1$, $\alpha<h(\beta)<\beta$.
Since $\alpha\leqone\beta$, this implies $\alpha\leqone h(\beta)$.
\qed

\blankline

The following lemma will be useful in establishing relations of the form $\alpha\leqtwoup\beta$.

\begin{lem}
\label{lblC1}
Assume $X< Y<\alpha \ \leqone \ \beta$ where $X$ and $Y$ are finite closed sets and let $X'$ be the set of $\xi\in X$ such that there is $\eta \in Y$ such that $\xi\leqtwo\eta$ .
If there are cofinally many $Y'$ below $\alpha$ such that $X\cup Y'$ is a covering of $X\cup Y$ and $\tilde{Y}$ is a subset of $\beta$ such that $X<\tilde{Y}$ and $X'\cup \tilde{Y}$ is a covering of $X'\cup Y$ then $X\cup \tilde{Y}$ is a covering of $X\cup Y$.
\end{lem}
{\bf Proof.}
Assume $\xi\in X$, $\eta \in Y$ and $\xi\leqone\eta$. 
We claim $\xi\leqone\beta$.
The lemma is easily verified using this fact.

Since there are cofinally many $Y'$ below $\alpha$ such that $X\cup Y'$ is a covering of $X\cup Y$, there are cofinally many $\eta'$ below $\alpha$ with $\xi\leqone\eta'$.
By part 2 of Lemma \ref{lblE}, $\xi\leqone\alpha$.
Since $\alpha\leqone\beta$, $\xi\leqone\beta$
\qed

\section{The Main Structural Lemma}
\label{smsl}

Recall that $\rho$ is an arbitrary additively indecomposable ordinal.

The following lemma is the key to the analysis in the remainder of the paper.

\blankline


\noindent{\bf Main Structural Lemma.}\
{\it Assume $\alpha\in ORD$ is not divisible by $\rho$. For all $\gamma\in ORD$, if $[0,\alpha)\not \leq_2 [\alpha,\alpha+ \gamma)$  then $ [1,1+\gamma]\cong   [\alpha,\alpha+\gamma] $.}

\blankline

\noindent{\bf Proof.}
Let $h$ be the operation $1+\eta\mapsto \alpha+\eta$. By Lemma \ref{lvB1},  $h$ is an isomorphism of $[1,\infty)$ and $[\alpha,\infty)$ as substructures of  \calRrho.

We will argue by induction on $\gamma$. Assume $\gamma\in ORD$ and the lemma is true for $\gamma'$ whenever $\gamma'<\gamma$.

Assume $[0,\alpha)\not\leq_2[\alpha,\alpha+\gamma)$. 
By the induction hypothesis, the restriction of $h$ is an isomorphism of $[1+\delta,1+\gamma)$ and $[\alpha,\alpha+\gamma)$.  
Since $[1,1+\gamma]\cong[\alpha,\alpha+\gamma]$ as substructures of \calRrho, it suffices to show that for $k=1,2$
 $$(\star) \ \ \ \ \ \ \ \ \ \ \ \   1+\delta \leq_k 1+\gamma \ \ \ {\rm iff} \ \ \ \alpha+\delta\leq_k \alpha+\gamma$$
whenever $\delta<\gamma$.

\halfblankline

{\bf Case 1:} Assume $k=1$. 

Suppose $\delta<\gamma$. We will show $(\star)$ holds. 

($\Rightarrow$) Assume $1+\delta\leqone 1+\gamma$. To show $\alpha+\delta \leqone \alpha+\gamma$, assume $X\subseteq \alpha+\delta$ and $Y\subseteq [\alpha+\delta,\alpha+\gamma)$ are finite with $X$ and $Y$ closed.  
By Lemma \ref{lblD1}, $1+\delta$ is divisible by $\rho$ implying $\alpha+\delta$ is divisible by $\rho$ (using the assumption that $\rho$ is additively indecomposable).
By Lemma \ref{lvA}, both $X$ and $Y$ are closed.
We will show there exists $\tilde{Y}\subseteq \alpha$ such that $X<\tilde{Y}$ and $X\cup\tilde{Y}$ is a covering of $X\cup Y$.
Let $X'$ be the collection of $\xi\in X$ such that $\xi\leqtwo\eta$ for some $\eta\in Y$.
By Lemma \ref{lblC}, it suffices to find $\tilde{Y}\subseteq \alpha$ such that $X<\tilde{Y}$ and $X'\cup\tilde{Y}$ is a covering of $X'\cup Y$.
Since $[0,\alpha)\not\leq_2[\alpha,\alpha+\gamma)$, $X'\subseteq [\alpha,\alpha+\delta)$.
Let $\mu$ be the largest element of $X\cup\{\alpha\}$.
Let $X'^*=h^{-1}[X']$, $Y^*=h^{-1}[Y]$  and $\mu^*=h^{-1}(\mu)$.
Since $\delta\leqone\gamma$, part 1 of Lemma \ref{lblA} implies there exists $\tilde{Y}^*\subseteq \delta$ such that $X'^*<\tilde{Y}^*$, $X'^*\cup\tilde{Y}^*$ is a covering of $X'^*\cup Y^*$ and $\mu^*<\tilde{Y}^*$. 
Let $\tilde{Y}=h[\tilde{Y}^*]$. 
Since $\mu^*<\tilde{Y}^*$, $\mu<\tilde{Y}$ implying $X<\tilde{Y}$.
Since the restriction of $h$ is an isomorphism of $[1,1+\gamma)$ and $[\alpha,\alpha+\gamma)$,  $X'\cup \tilde{Y}$ is a covering of $X'\cup Y$. 

($\Leftarrow$) The proof is similar, but easier, than the $(\Rightarrow)$ direction.

\halfblankline

{\bf Case 2:} Assume $k=2$.

Suppose $\delta<\gamma$.

($\Rightarrow$) Assume $1+\delta\leqtwo 1+\gamma$.

By Lemma \ref{lblD1}, $1+\delta$ is of the form $\rho\cdot\lambda$ for a limit ordinal $\lambda$. 
This implies $\delta=\rho\cdot\lambda$.

By the case $k=1$, $\alpha+\delta\leqone \alpha+\gamma$.

We omit the proof that $\rho\cdot\alpha+\delta\leqtwodown \rho\cdot\alpha+\gamma$ since it is similar to  the proof of Case 1 (using the fact that we have established the case $k=1$ to handle the additonal condition in the definition of $\leqtwo$).

To show that $\alpha+\delta \leqtwoup \alpha+\gamma$, assume that $X$ and $Y$ are finite  subsets of $\alpha+\delta$ with $X$ and $Y$ closed and $X<Y$ such that there are cofinally many $\tilde{Y}$ below $\alpha+\delta$ such that $X\cup \tilde{Y}$ is a covering of $X\cup Y$. 
We will show that there are cofinally many $\tilde{Y}$ below $\alpha+\gamma$ such that $X\cup \tilde{Y}$ is a covering of $X\cup Y$. 
Let $X'$ be the collection of $\xi\in X$ such that $\xi\leqtwo\eta$ for some $\eta\in Y$.
Since $[0,\alpha)\not\leq_2 [\alpha,\alpha+\gamma)$, $X'\subseteq [\alpha,\alpha+\delta)$.
By Lemma \ref{lblC1}, it suffices to show there are cofinally many $\tilde{Y}$ below $\alpha+\gamma$ such that $X'\cup \tilde{Y}$ is a covering of $X'\cup Y$. 
This is straightforward using the induction hypothesis.

$(\Leftarrow)$ The proof is similar to the $(\Rightarrow)$ direction but easier.
\qed

\begin{lem}
\label{lmslA}
Assume $\alpha_1$ and $\alpha_2$ are are ordinals and $\gamma\in ORD$ satisfy $[0,\rho\cdot\alpha_i] \not \leq_2 [\rho\cdot\alpha_i+1,\rho\cdot\alpha+ \gamma)$ for $i=1,2$.
Moreover, let $\delta_i$ satisfy $\rho\cdot\alpha_i+\delta_i=max_1(\rho\cdot\alpha_i)$ for $i=1,2$ ($\delta_i$ may be $\infty$).
\begin{enumerate}
\item
If  either $\gamma\leq\delta_2$ or $\delta_1\leq\delta_2$ then
$[\rho\cdot\alpha_2,\rho\cdot\alpha_2+ \gamma]$ is a covering of $[\rho\cdot\alpha_1,\rho\cdot\alpha_1+ \gamma]$.
\item
If  either $\gamma\leq\delta_1,\delta_2$ or $\delta_1=\delta_2$ then 
$$[\rho\cdot\alpha_1,\rho\cdot\alpha_
1+ \gamma]\cong[\rho\cdot\alpha_2,\rho\cdot\alpha_2+ \gamma]$$
\end{enumerate}
\end{lem}
{\bf Proof.}
By Lemma \ref{lvB} 
$$[\rho\cdot\alpha_1,\rho\cdot\alpha_
1+ \gamma]\cong[\rho\cdot\alpha_2,\rho\cdot\alpha_2+ \gamma]$$
 as substructures of \calRrho, and by the Main Structural Lemma  
$$[\rho\cdot\alpha_1+1,\rho\cdot\alpha_
1+ \gamma]\cong[\rho\cdot\alpha_2+1,\rho\cdot\alpha_2+ \gamma]$$
 as substructures of \calRtworho. 
 To verify part 1, it remains to show that for $k=1,2$ and all $\xi\leq\gamma$, if $\rho\cdot\alpha_1\leq_k\rho\cdot\alpha_1+\xi$ then $\rho\cdot\alpha_2\leq_k\rho\cdot\alpha_2+\xi$.
By assumption, the case $k=2$ is vacuous. The case $k=1$ follows from the assumption $\gamma\leq\delta_2$ or $\delta_1\leq\delta_2$.

Part 2 follows from part 1.
 \qed

\section{The First Recurrence Theorem for \leqone}
\label{sfirstrone}

Recall that $\rho$ is an arbitrary additively indecomposable ordinal.

The decomposition of $ORD$ into connectivity components with respect to \leqone\ provides insight into the structure of \calRtworho.

\begin{dfn}
 Assume $\alpha\in ORD$ and $k\in\{1,2\}$. If there is a bound on the ordinals $\beta$ such that $\alpha\leqk\beta$ then $max_k^\rho(\alpha)$ is the largest $\beta \in ORD$ such that $\alpha \leq_k \beta$. Otherwise, we write $max_k^\rho(\alpha)=\infty$.
\end{dfn}

Since $\rho$ is fixed we will often omit $\rho$ as a parameter in notation as in the previous section.


\begin{lem}
\label{lgsA}
The collection of $\kappa\in ORD$ which are minimal with respect to \leqone\ is topologically closed.
\end{lem}
{\bf Proof.}
Assume $\kappa\in ORD$ is a limit of minimal elements of $ORD$ with respect to \leqone. To see that $\kappa$ is minimal with respect to \leqone, argue by contradiction and assume $\xi <_1 \kappa$. There is some $\eta$ such that $\xi < \eta < \kappa$ and $\eta$ is minimal with respect to \leqone. Since \leqone\ respects $\leq$, $\xi <_1 \eta$ contradicting the fact that $\eta$ is minimal with respect to $\leqone$. \qed

\blankline

The lemma implies that if the ordinals which are minimal with respect to \leqone\ are bounded then there is a largest such ordinal.

\begin{dfn}
Let $\kappa_\alpha^\rho$ $(\alpha \in \theta_1^\rho)$ enumerate  the elements of $ORD$  which are minimal with respect to \leqone\  (we allow the possibility that $\theta_1^\rho = \infty$). Define  $I_\alpha^\rho$ to be the collection of $\beta$ such that $\kappa_\alpha\leqone\beta$ when $\alpha \in \theta_1^\rho$.
\end{dfn}

Notice that  $I_\alpha=[\kappa_\alpha,max_1(\kappa_\alpha)]$  when $max_1(\kappa_\alpha)\not=\infty$, and $I_\alpha=[\kappa_\alpha,\infty)$ otherwise.

In mildly strong theories, ${\sf ZF}$ is much stronger than required, there is an ordinal $\kappa$ with $\kappa \leqone \infty$. The least such $\kappa$ is the largest ordinal which is minimal with respect to  $\leqone$. 


\begin{lem}
\
\label{lgsB}
\begin{enumerate}
 \item 
   $\kappa_0=0$ and $max_1(\kappa_0)=\kappa_0$. 
 \item 
 $\alpha \mapsto \kappa_\alpha$ is continuous.
 \item 
 If $\theta_1\not=\infty$ then there is an ordinal $\theta$ such that $\theta+1=\theta_1$ and  $max_1(\kappa_\theta)=\infty$. 
 \item
 For all $\alpha<\theta_1$, $[0,\kappa_\alpha)\not\leq_1 [\kappa_\alpha,\infty)$.
 \item 
 If $\alpha < \beta < \theta_1$ then $I_\alpha < I_\beta$ and $I_\alpha\not\leq_1 I_\beta$.
  \item 
  $ORD = \bigcup _{\xi \in\theta_1} I_\xi$.
 \item
  If $\alpha \in\theta_1$ then $\kappa_\alpha = \bigcup_{\xi<\alpha}I_\xi$.
  \item
   If $\alpha+1<\theta_1$ and $0<\xi<\rho\cdot\omega$ then $\alpha+\xi<\theta_1$, $\kappa_{\alpha+\xi}=max_1(\kappa_\alpha) + \xi$ and $max_1(\kappa_{\alpha+\xi})=\kappa_{\alpha+\xi}$.
\item
If $\alpha+1<\theta_1$ then $\alpha+\rho\cdot\omega<\theta_1$ and $\kappa_{\alpha+\rho\cdot\omega}=\kappa_\alpha+\rho\cdot\omega$.
\item
If $\alpha<\theta_1$ then $\kappa_\alpha$ is divisible by $\rho$ iff $\alpha$ is divisible by $\rho$.
 \item 
 For all $\alpha<\theta_1$, there is no $\beta$ such that $\kappa_\alpha <_2 \beta$.
 \end{enumerate}
\end{lem}
{\bf Proof.}
Parts 1-7 follow directly from the definitions and Lemma \ref{lgsA}.

For part 8, assume $\alpha+1<\theta_1$.  Lemma \ref{lblD1} implies that $max_1(\kappa_\alpha+\xi)=\kappa_\alpha+\xi$ whenever $0<\xi<\rho\cdot\omega$. This implies $\kappa_{\alpha+\xi}=max_1(\kappa_\alpha)+\xi$ by induction on $\xi<\rho\cdot\omega$.

Part 9 follows from parts 2 and 8.

Part 10 follows by induction on $\alpha$ using parts 1, 2 and 8 (notice that $\beta+\rho$ is always divisible by $\rho$ since $\rho$ is additively indecomposable).

Part 11 follows from part 1 of Lemma \ref{lblD}.
\qed

\blankline

We remark that  $max_1(\rho\cdot\omega)=\rho\cdot\omega +1$.
A more general result will be proved later.

In the two recurrence theorems for \leqone, we will determine when $I_\alpha\cong I_\beta$. We can make some simple observations here. 

By part 8 of the lemma, $\kappa_{\alpha+\xi}$ is not divisible by $\rho$ when $0<\xi<\rho$ and $I_{\alpha+\xi}=\{\kappa_{\alpha+\xi}\}$ which is not a closed set. This easily implies that the intervals $I_{\alpha+\xi}$ with $0<\xi<\rho$ are all isomorphic to each other.

On the other hand, $\kappa_{\rho\cdot\beta}$ is always divisible by $\rho$ by part 10 of the lemma.
This implies $I_{\rho\cdot\beta}$ is closed.
Therefore, $I_{\rho\cdot\beta}$ is never isomorphic to $I_{\alpha+\xi}$ when $0<\xi<\rho$.
We are left with determining when $I_{\rho\cdot\beta_1}\cong I_{\rho\cdot\beta_2}$.

The following special case of Lemma \ref{lmslA} will be particularly useful.

\begin{lem}
\label{lgsB1}
For all $\alpha$ with $\alpha+1<\theta_1$, $[\kappa_{\alpha+\rho},\infty)\cong \calRtworho$.
\end{lem}
{\bf Proof.} By parts 10 and 8 of the previous lemma, $\kappa_{\alpha+\rho}$ is divisible by $\rho$ and $max_1(\kappa_{\alpha+\rho})=\kappa_{\alpha+\rho}$.
In particular, $[0,\kappa_{\alpha+\rho}]\not\leq_2 [\kappa_{\alpha+\rho}+1,\infty)$.
Since $0$ is divisible by $\rho$ and $max_1(0)=0$, part 2 of Lemma \ref{lmslA} implies the desired conclusion.
\qed


\begin{lem}
\label{lgsD}
Assume $\alpha+1<\theta_1$ and $0<\beta<\theta_1$.
\begin{enumerate}
\item 
$\alpha+\beta<\theta_1$.
\item 
$\kappa_{\alpha+\beta}=max_1(\kappa_\alpha)+\kappa_\beta$
\item
$max_1(\kappa_{\alpha+\beta})=max_1(\kappa_\alpha)+max_1(\kappa_\beta)$
\end{enumerate}
\end{lem}
{\bf Proof.} 
By part 8 of Lemma \ref{lgsB}, $\kappa_{\alpha+1}=max_1(\kappa_\alpha)+1$. 
By the Main Structural Lemma, the operation $h$ given by $\gamma\mapsto max_1(\kappa_\alpha)+\gamma$ is an isomorphism of $[1,\infty)$ and $[\kappa_{\alpha+1},\infty)$. 
An ordinal in $[1,\infty)$ is a minimal element of $[1,\infty)$ with respect to \leqone\ iff it is a minimal element of $ORD$ with respect to $\leqone$. 
The collection of such ordinals is enumerated by $\kappa_{1+\xi}$ $(1+\xi\in\theta_1)$. 
Similarly, an ordinal in $[\kappa_{\alpha+1},\infty)$ is a minimal element of $[\kappa_{\alpha+1},\infty)$ with respect to \leqone\ iff it is a minimal element of $ORD$ with respect to \leqone.
The collection of such ordinals is enumerated by $\kappa_{\alpha+1+\xi}$ $(\alpha+1+\xi\in\theta_1)$.
Since $h$ is an isomoprhism, $h$ maps the ordinals which are minimal in $[1,\infty)$ with respect to \leqone\ onto the ordinals which are minimal in $[\kappa_{\alpha+1},\infty]$ with respect to $\leqone$. 
Therefore, $h(\kappa_{1+\xi})=\kappa_{\alpha+1+\xi}$ whenever $1+\xi\in\theta_1$. 

To establish parts 1 and 2, choose $\xi$ so that $\beta=1+\xi$. 
By the previous paragraph, $h(\kappa_\beta)=\kappa_{\alpha+\beta}$.
Using the definition of $h$, $max_1(\kappa_\alpha)+\kappa_{\beta}=\kappa_{\alpha+\beta}$. 

Using the fact $h$ is an isomorphism again, $h(max_1(\kappa_\beta))=max_1(h(\kappa_\beta))$. 
Using the definition of $h$ and part 2, $max_1(\kappa_\alpha)+max_1(\kappa_\beta)=max_1(\kappa_{\alpha+\beta})$.
\qed

\begin{lem}
\label{lgsDD}
Either $\theta_1=\infty$ or $\theta_1=\theta+1$ where $\theta$ is additively indecomposable and greater than $\rho$.
\end{lem}
{\bf Proof.}
Assume $\theta_1\not=\infty$.
By part 3 of Lemma \ref{lgsB}, there exists $\theta$ such that $\theta_1=\theta+1$.
By part 9 of Lemma \ref{lgsB}, $\theta=\rho\cdot\lambda$ for some limit ordinal $\lambda$.
By the previous lemma, if $\alpha<\theta$ then $\alpha+\theta=\theta$.
Therefore, $\theta$ is additively indecomposable. 
\qed

\blankline

We will see later that $\theta_1$ is either $\infty$ or the successor of an epsilon number greater than $\rho$.

\begin{lem}
\label{lgsD0}
Assume $\alpha<\theta_1$.
$\kappa_\alpha$ is additively indecomposable iff $\alpha$ is additively indecomposable.
\end{lem}
{\bf Proof.}
Follows easily from part 2 of Lemma \ref{lgsD} and the continuity of $\alpha\mapsto \kappa_\alpha$ (part 2 of Lemma \ref{lgsB}).
\qed

\begin{thm}
\label{tgsD1}
{\bf (First Recurrence Theorem for \leqone)}
If $\alpha+1<\theta_1$ and $0<\beta<\theta_1$ then $I_{\alpha+\beta}\cong I_\beta$.
\end{thm}
{\bf Proof.}
Assume $\alpha+1<\theta_1$ and $0<\beta<\theta_1$.

By the Main Structural Lemma, the map $h$ given by $h(\xi)=max_1(\kappa_\alpha)+\xi$ is an isomorphism of $[1,\infty)$ and $[max_1(\kappa_\alpha)+1,\infty)$. 
Lemma \ref{lgsD} implies
$$h(\kappa_\beta)=max_1(\kappa_\alpha)+\kappa_\beta=\kappa_{\alpha+\beta}$$

First, assume $\beta+1<\theta_1$.  $I_\beta=[\kappa_\beta,max_1(\kappa_\beta)]\cong [h(\kappa_\beta),h(max_1(\kappa_\beta))]$.
By Lemma \ref{lgsD} again,
 $$h(max_1(\kappa_\beta))=max_1(\kappa_\alpha)+max_1(\kappa_\beta)=max_1(\kappa_{\alpha+\beta})$$
 Therefore, $I_\beta\cong I_{\alpha+\beta}$.
 
 Now, assume $\beta+1=\theta_1$. 
By Lemma \ref{lgsDD}, $\beta$ is additively indecomposable and greater than $\rho$. 
Since $\alpha<\beta$, $\alpha+\beta=\beta$ making the conclusion of the theorem trivial.
\qed

\section{Incompressible Sets}
\label{sis}

Recall that $\rho$ is an arbitrary additively indecomposable ordinal.

In this section, we introduce a weakening of the notion of isominimal set from [\ref{Ca09}] sufficient for this paper. 
The existence of these {\it incompressible sets} is easier to establish than that of isominimal sets and the inclusion of this section makes the paper self-contained.
On the other hand, for those familiar with [\ref{Ca09}], this section can be skipped and the notion of an incompressible covering of a set $Y$ can be replaced by the notion of an isominimal copy of $Y$ and the notion of an incompressible set can be replaced by the notion of an isominimal set.

\begin{lem}
\label{lgsDA}
Assume $\rho\cdot\delta< \theta_1$ and $\delta$ is an infinite additively indecomposable ordinal.
If $Y$ is a finite closed subset of $\kappa_{\rho\cdot\delta}$ then there are cofinally many subsets of $\kappa_{\rho\cdot\delta}$ which are isomorphic to $Y$.
\end{lem}
{\bf Proof.}
Assume $Y$ is a finite closed subset of $\kappa_{\rho\cdot\delta}$. 
To show there are cofinally many subsets of $\kappa_{\rho\cdot\delta}$ isomorphic to $Y$, suppose $\xi<\kappa_{\rho\cdot\delta}$. 
There exists $\gamma<\delta$ such that $Y$ is a subset of $\kappa_{\rho\cdot\gamma}$ and $\xi<\kappa_{\rho\cdot\gamma}$.
By  part 2 of Lemma \ref{lmslA} and part 2 of Lemma  \ref{lgsD}, $[\kappa_{\rho\cdot(\gamma+1)},\kappa_{\rho\cdot(\gamma+1)+\rho\cdot\gamma})$ is isomorphic to $[0,\kappa_{\rho\cdot\gamma})$. 
Therefore, there is a subset $\tilde{Y}$ of $[\kappa_{\rho\cdot(\gamma+1)},\kappa_{\rho\cdot(\gamma+1)+\rho\cdot\gamma})$ which is isomorphic to $Y$. 
By choice of $\gamma$, $\xi<\tilde{Y}$.
By assumption, $\delta$ is closed under addition. 
Therefore, $\gamma+1+\gamma<\delta$ implying $\tilde{Y}\subseteq \kappa_{\rho\cdot\delta}$. 
\qed

\begin{dfn}
For an ordinal $\beta$, define the {\bf index} of $\beta$, $index^{\rho}(\beta)$ to be the unique $\alpha<\theta_1$ such that $\beta\in I_\alpha$.
\end{dfn}

As usual, we will write $index(\beta)$ for $index^\rho(\beta)$ when $\rho$ is clear from the context.

\begin{dfn}
Assume $Y$ is a finite closed set of ordinals. 
A covering $h$ of $Y$ in $\calRtworho$ is a {\bf $\rho$-incompressible covering} of $Y$ if $index(h(\beta))\leq index(h'(\beta))$ for all $\beta\in Y$ whenever $h'$ is a covering of $Y$ in $\calRtworho$.
The range of a $\rho$-incompressible covering of $Y$ will also be called a $\rho$-incompressible covering of $Y$.
$Y$ is {\bf $\rho$-incompressible} if $Y$ is a $\rho$-incompressible covering of itself.
\end{dfn}

As usual, we will drop mention of $\rho$ when it is understood from context and write ``incompressible covering" and ``incompressible" for ``$\rho$-incompressible covering" and ``$\rho$-incompressible" respectively.

Notice that being an incompressible covering of $Y$ is not the same as being a covering of $Y$ which is incompressible.

\begin{lem}
\label{lisA}
\begin{enumerate}
\item
Any finite union of incompressible sets is incompressible.
\item
If $h$ and $h'$ are incompressible coverings of $X$ then $index(h(\beta))=index(h'(\beta))$ for all $\beta\in X$.
\item
If $X$ is incompressible and $X\cap I_\alpha\not=\emptyset$ then $X\cup \{\kappa_\alpha\}$ is incompressible.
\item
Assume $X$ is incompressible and $\kappa_\alpha\in X$ whenever $X\cap I_\alpha\not=\emptyset$. If $Y$ is a finite closed set and $index[Y]\subseteq index[X]$ then $X\cup Y$ is incompressible.
\item
Assume $X$ is a nonempty finite closed set such that $min(X)\leqone max(X)$.
A covering $\tilde{X}$ of $X$ is an incompressible covering of $X$ iff $\tilde{X}\subseteq I_\alpha$ where $\alpha$ is minimal such that $I_\alpha$ contains a covering of $X$.
\item
Assume $X$ and $Y$ are finite  sets such that $X\cup Y$ is closed, no element of $Y$ is indecomposable, $X<Y$ and $X\not\leq_1 Y$. If $h$ is a covering of $X\cup Y$ then $h$ is an incompressible covering of $X\cup Y$ iff $h\!\!\upharpoonright\!\!X$ is an incompressible covering of $X$.
\item
Assume $X$ and $Y$ are finite nonempty closed sets such that $X<Y$ and $X\not\leq_1 Y$. If $\tilde{X}$ is a covering of $X$ and $\tilde{Y}$ is a covering of $Y$ such that $\tilde{X}<\tilde{Y}$ then $\tilde{X}\cup\tilde{Y}$ is an incompressible covering of $X\cup Y$ iff the following conditions hold:
\begin{enumerate}
\item
$\tilde{X}$ is an incompressible covering of $X$.
\item
$\tilde{Y}=max(\tilde{X})+\rho+Y^*$ where $Y^*$ is an incompressible covering of $Y$.
\end{enumerate}
\item
Assume $X$ and $Y$ are finite sets  such that $X\cup Y$ is closed, $X<Y$ and $X\not \leq_1 Y$.
If $h$ is an incompressible covering of $X\cup Y$  then the restriction of $h$ to $X$ is an incompressible covering of $X$. In particular, if $X\cup Y$ is incompressible then so is $X$.
\item
If $h$ is an incompressible covering of $X$  and $h[X]\cap I_\alpha\not=\emptyset$ then $\alpha+1<\theta_1$, $max(I_\alpha)\in h[X]$ and if $h(\beta)=max(I_\alpha)$ then $h(\beta)\leq h'(\beta)$ for any covering $h'$ of $X$.
In particular, if $X$ is incompressible and $X\cap I_\alpha\not=\emptyset$ then $\alpha+1<\theta_1$, $max(I_\alpha)\in X$ and if $\beta=max(I_\alpha)$ then $\beta\leq h'(\beta)$ for any covering $h'$ of $X$.
\end{enumerate}
\end{lem}
{\bf Proof.} 
Parts 1 and 2 are immediate.

For part 3, assume $h$ is a covering of $X\cup \{\kappa_\alpha\}$. Let $\beta\in I_\alpha$. 
Since $X$ is incompressible, $h(\beta)\in I_\xi$ where $\alpha\leq \xi$.
Since $\kappa_\alpha\leqone \beta$, $h(\kappa_\alpha)\leqone h(\beta)$ implying $h(\kappa_\alpha)\in I_\xi$.

Parts 4 is straightforward. 

Part 5 is straightforward after noticing that any covering of $X$ must be contained in some $I_\xi$.

For parts 6 and 7 we will use

\halfblankline

{\bf Claim.} 
\vspace{-4mm}
 \begin{enumerate}
\item
Assume $h$ is a covering of $X$ and $\beta=max(X)$.
If $h(\beta)\in I_\alpha$ and $h(\beta)$ is not the largest element of $I_\alpha$ then there is a covering $h'$ of $X$ such that $h'\leq h$ and $h':X\rightarrow \kappa_\alpha$.
\item
If $h$ is an incompressible covering of $X$, $\beta=max(X)$ and $h(\beta)\in I_\alpha$ then $h(\beta)=max(I_\alpha)$ and $h(\beta)\leq h'(\beta)$ for any covering $h'$ of $X$.
\item If $h$ is an incompressible covering of $X$, $\beta$ is the largest indecomposable in $X$ and $h(\beta)\in I_\alpha$ then $h(\beta)$ is the largest indecomposable in $I_\alpha$ and $h(\beta)\leq h'(\beta)$ for every covering $h'$ of $X$.
\end{enumerate}

Part 1 of the claim is immediate from the definition of $\leqone$. Parts 2 and 3 follow from part 1.

\halfblankline

For part 6 of the lemma, suppose $h$ is a covering of $X\cup Y$ and let $h^-$ be the restriction of $h$ to $X$.

First, assume $h$ is an incompressible covering of $X\cup Y$. 
To show $h^-$ is an incompressible covering of $X$, assume $f$ is a covering of $X$.
By Lemma \ref{lvBB}, every embedding of $X$, as a substructure of \calRrho, into \calRrho\ extends uniquely to an embedding of $X\cup Y$.
Let $h'$ be the embedding of $X\cup Y$, as a substructure of \calRrho, into \calRrho\ which extends $f$. 
Clearly, $h'$ is a covering of $X\cup Y$.
Since $h$ is incompressible, $index(h(\alpha))\leq index(h'(\alpha))$ for all $\alpha$ in $X\cup Y$.
In particular, $index(h^-(\alpha))\leq index(f(\alpha))$ for all $\alpha\in X$.

Now assume that the restriction of $h$ to $X$ is an incompressible covering of $X$.
To show $h$ is an incompressible covering of $X\cup Y$, assume $h'$ is a covering of $X\cup Y$.
By assumption, $index(h(\alpha))\leq index(h'(\alpha))$ for all $\alpha\in X$.
Let $\beta$ be the largest indecomposable ordinal in $X$.
By part 3 of the of the claim above, $h(\beta)\leq h'(\beta)$.
Suppose $\alpha\in Y$. 
Since $X\cup Y$ is closed, there exists $\xi<\rho$ such that $\alpha=\beta+\xi$.
Therefore, $h(\alpha)=h(\beta)+\xi\leq h'(\beta)+\xi=h'(\alpha)$ implying $index(h(\alpha))\leq index(h'(\alpha))$.

\halfblankline

For part 7, assume $\tilde{X}$ is a covering of $X$ and $\tilde{Y}$ is a covering of $Y$ such that $\tilde{X}<\tilde{Y}$. 
Let $\sigma=index(max(\tilde{X}))$.  

\halfblankline

{\bf Claim for part 7.} Assume $\tilde{X}$ is an incompressible covering of $X$. For all $\beta$, $index(max(\tilde{X})+\rho +\beta)=\sigma+\rho+index(\beta)$.

By part 2 of the claim above, $max(\tilde{X})=max(I_\sigma)=max_1(\kappa_\sigma)$.
By part 8 of Lemma \ref{lgsB}, $max(\tilde{X})+\rho=\kappa_{\sigma+\rho}$ and $max_1(\kappa_{\sigma+\rho})=\kappa_{\sigma+\rho}$.
By Lemma \ref{lgsD}, $\kappa_{\sigma+\rho+\xi}=\kappa_{\sigma+\rho}+\kappa_\xi=max(\tilde{X})+\rho+\kappa_\xi$ for all $\xi$ with $0<\xi<\theta_1$.
This equation also holds for $\xi=0$.
Therefore, $\kappa_{\sigma+\rho+\xi}\leq max(\tilde{X})+\rho+\beta$ iff $\kappa_\xi \leq \beta$ for all $\xi<\theta_1$.
This implies the conclusion of the claim.

\halfblankline

$(\Rightarrow)$ Suppose $\tilde{X}\cup \tilde{Y}$ is an incompressible covering of $X\cup Y$.
Let $h$ be the covering of $X\cup Y$ with range $\tilde{X}\cup \tilde{Y}$.

To verify (a), suppose $f$ is a covering of $X$. 
We will show $index(h(\beta))\leq index(f(\beta))$ for all $\beta\in X$. 
By part 1 of the claim above, we may assume the range of $f$ is contained in $\kappa_\alpha$ for some $\alpha<\theta_1$.
By increasing $\alpha$ if necessary, we may assume $\tilde{X}\cup \tilde{Y}<\kappa_\alpha$ by part 2 of the claim above. 
Finally, we may assume $\alpha$ has the form $\rho\cdot\delta$ where $\delta$ is infinite and additively indecomposable by Lemma \ref{lgsDD}.
Let $X'$ be the range of $f$.
By Lemma \ref{lgsDA}, there exists $Y'\subseteq \kappa_\alpha$ which is isomorphic to $\tilde{Y}$ such that $X'<Y'$. 
Clearly, $X'\cup Y'$ is a covering of $X\cup Y$ and the covering $h'$ of $X\cup Y$ with range $X'\cup Y'$ extends $f$.
Using the fact $h$ is incompressible, $index(h(\beta))\leq index(h'(\beta))=index(f(\beta))$  for $\beta\in X$.

To verify (b), notice that the least element of $\tilde{Y}$ must be indecomposable and greater than the largest element of $\tilde{X}$.
Therefore, $max(\tilde{X})+\rho\leq min(\tilde{Y})$ and there is a finite closed set $Y^*$ such that $\tilde{Y}=max(\tilde{X})+\rho+Y^*$.

To show that $Y^*$ is incompressible, assume $Y'$ is a covering of $Y$.
Suppose $i<card(Y)$. 
We need to show $index((Y^*)_i)\leq index((Y')_i)$.
By Lemma \ref{lgsB1}, $max(\tilde{X})+\rho+Y'$ is isomorphic to $Y'$ implying it is a covering of $Y$.
This implies $\tilde{X}\cup (max(\tilde{X})+\rho + Y')$ is a covering of $X\cup Y$.
Since $\tilde{X}\cup \tilde{Y}$ is an incompressible covering of $X\cup Y$, $index((max(\tilde{X})+\rho+ Y^*)_i)\leq index((max(\tilde{X})+rho+Y')_i)$.
By the claim for part 7, $\sigma+\rho + index((Y^*)_i)\leq \sigma+\rho+ index((Y')_i)$. 
Therefore, $index((Y^*)_i)\leq index((Y')_i)$.

$(\Leftarrow)$ Assume (a) and (b). 
To show $\tilde{X}\cup \tilde{Y}$ is an incompressible covering of $X\cup Y$, assume $\overline{X}\cup \overline{Y}$ is a covering of $X\cup Y$ where $\overline{X}$ is a covering of $X$, $\overline{Y}$ is a covering of $Y$ and $\overline{X}<\overline{Y}$.
Since $\tilde{X}$ is an incompressible covering of $X$, $index((\tilde{X})_i)\leq index((\overline{X})_i)$ for all $i<card(X)$.
It remains to show $index((\tilde{Y})_i)\leq index((\overline{Y})_i)$ for all $i<card(Y)$.

Fix $i<card(Y)$
Since the least element of $\overline{Y}$ is indecomposable, $max(\overline{X})+\rho\leq \overline{Y}$.
By part 2 of the claim above, $max(\tilde{X})\leq max(\overline{X})$.
Therefore, $max(\tilde{X})+\rho\leq \overline{Y}$ implying $\overline{Y}=max(\tilde{X})+\rho +Y'$ for some closed $Y'$.
By Lemma \ref{lgsB1}, $Y'$ is isomorphic to $\overline{Y}$ implying it is a covering of $Y$.
Since $Y^*$ is an incompressible covering of $Y$, $index((Y^*)_i)\leq index((Y')_i)$.
By the claim for part 7, $index((\tilde{Y})_i)\leq index((\overline{Y})_i)$.

Part 8 follows from parts 6 and 7.

For part 9, let $X_0$ be the collection of $\beta\in X$ such that $h(\beta)\leq max(I_\alpha)$.
Part 8 implies the restriction of $h$ to $X_0$ is an incompressible covering of $X_0$.
The conclusion of part 9 follows from part 2 of the claim.
\qed

\begin{lem}
\label{lisB}
Assume $X$ is a finite closed set of ordinals.
There is an incompressible covering of $X$.
\end{lem}
{\bf Proof.}
By induction on the cardinality of $X$.
When $X$ is empty, the lemma is trivial.
So, we may assume $X$ is nonempty.
Let $X_1$ be the collection of $\beta\in X$ such that $min(X)\leqone \beta$.
Let $X_2$ be the collection of $\beta\in X$ such that $max(X_1)<\beta<max(X_1)+\rho$.
Let $X_3$ be the collection of $\beta\in X$ such that $max(X_1)+\rho\leq\beta$.
Clearly, $X_1<X_2<X_3$, $X_1\not\leq_1 X_2\cup X_3$ and $X_2\not\leq_1 X_3$.
If $X_3\not=\emptyset$, the existence of an incompressible covering follows from the induction hypothesis and part 7 of the previous lemma. If $X_3=\emptyset$ and $X_2\not=\emptyset$, the existence of an incompressible covering follows from part 6 of the previous lemma. 
If $X_2=X_3=\emptyset$, the existence of an incompressible covering follows from part 5 of the previous lemma.
\qed

\begin{lem}
\label{lgsF}
If $K$ is a finite set of ordinals such that 
\begin{enumerate}
\item
If $\alpha\in K$ then $\alpha+1<\theta_1$.
\item
If $\rho\cdot\delta+\xi\in K$ and $\xi<\rho$ then $\rho\cdot\delta\in X$.
\item
If $\rho\cdot (n+1)\in K$ where $n\in\omega$ then 
$$\{0,\rho\cdot 1, \rho\cdot 2,\ldots,\rho\cdot n\}\subseteq K$$
\item
If $n\in\omega$ and $\rho\cdot(\delta_0+\cdots +\delta_{n+1})\in K$ where $\delta_0\geq\cdots\geq\delta_{n+1}$ are additively indecomposable and $1<\delta_0$ then 
$$\{\rho\cdot \delta_0,\rho\cdot(\delta_0+\delta_1),\ldots, \rho\cdot(\delta_0+\cdots+\delta_n)\}\subseteq K$$
\end{enumerate}
then there exists an incompressible set $X$ with $index[X]=K$.
\end{lem}
{\bf Proof.} 
We begin with a special case.

\halfblankline

{\bf Claim 1.} Assume $\rho\cdot\delta+1<\theta_1$ where $\delta$ is additively indecomposable and $1<\delta$. There is an incompressible $Y$ with $Y\subseteq I_{\rho\cdot\delta}$.

Choose finite closed $X\subseteq \kappa_{\rho\cdot\delta}$ and finite closed $Y\subseteq [\kappa_{\rho\cdot\delta},max_1(\kappa_{\rho\cdot\delta})]$ with the property that there is no $\tilde{Y}\subseteq \kappa_{\rho\cdot\delta}$ such that $X<\tilde{Y}$ and $X\cup \tilde{Y}$ is a covering of $X\cup Y$. 
We may assume $\kappa_{\rho\cdot\delta}\in Y$.
Argue by contradiction and assume $Y$ is not incompressible. 
Let $\tilde{Y}$ be a covering of $Y$ such that $\beta=index((\tilde{Y})_i)<index((Y)_i)$ for some $i$. 
Since $min(Y)=\kappa_{\rho\cdot\delta}\leqone max(Y)$, $min(\tilde{Y})\leqone max(\tilde{Y})$. 
Therefore, $\tilde{Y}\subseteq I_\beta$.
By the previous lemma, we may assume that $X<\tilde{Y}$.
Since $\xi\not\leq_1\eta$ for any $\xi\in X$ and $\eta \in Y$, $X\cup\tilde{Y}$ is a covering of $X\cup Y$ -- contradiction.

\halfblankline

Since any finite union of incompressible sets is incompressible by part 1 of Lemma \ref{lisA}, the following claims imply the theorem.

\halfblankline

{\bf Claim 2.} For $n\in \omega$ and $\xi<\rho$, $X=\{0,\rho\cdot 1,\rho\cdot 2,\ldots,\rho\cdot n,\rho\cdot n+\xi\}$ is incompressible and $index[X]=X$.

Straightforward since any covering maps indecomposable ordinals to indecomposable ordinals. Notice that $index(\rho\cdot i)=\rho\cdot i$  for all $i\in \omega$ and $index(\rho\cdot n +\xi)=\rho\cdot n +\xi$ by parts 1 and 8 of Lemma \ref{lgsB}. 

\halfblankline

{\bf Claim 3.}
Suppose $\delta_0\geq \delta_1\geq\cdots \geq \delta_n$ are additively indecomposable, $1<\delta_0$ and $\xi<\rho$.
There is an incompressible set $X$ with 
$$index[X]=\{\rho\cdot \delta_0,\rho\cdot(\delta_0+\delta_1),\ldots, \rho\cdot(\delta_0+\cdots+\delta_n),\rho\cdot(\delta_0+\cdots+\delta_n)+\xi\}$$

For $i\leq n$, let $\delta_i^-$ satisfy $\delta_i=1+\delta_i^-$. 
So, if $\delta_i=1$ then $\delta_i^-=0$ and $\delta_i^-=\delta_i$ otherwise. 
By Claim 1 and the fact that $\{0\}$ is an incompressible subset of $I_0$, there exists an incompressible subset  $X_i^*$  of $I_{\rho\cdot\delta_i^-}$.
For $i\leq n$, define $X_i$ by induction so that $X_0=X_0^*$ and $X_i=max(X_j)+\rho+X_i^*$ when $i=j+1$. 
Clearly, $X_0<X_1<\cdots<X_n$.

We claim $X_0\cup\cdots\cup X_i$ is incompressible and $X_i\subseteq I_{\rho\cdot(\delta_0+\cdots+\delta_i)}$ for $i=0,\ldots,n$.

Argue by induction. 

The case $i=0$ is clear by choice of $X_0^*$. 

Assume $i=j+1$. 
By part 7 of Lemma \ref{lisA}, $X_0\cup\cdots\cup X_i$ is incompressible. 
Since 
$$\kappa_{\delta_i^-}\leq X_i^*<\kappa_{\delta_i^-+1}$$
we have
 $$max(X_j)+\rho+\kappa_{\rho\cdot\delta_i^-}\leq X_i<\max(X_j)+\rho+\kappa_{\rho\cdot\delta_i^-+1}$$
By part 9 of Lemma \ref{lisA}, $max(X_j)=max_1(X_0\cup \cdots \cup X_j)=max_1(\kappa_{\rho\cdot(\delta_0+\cdots+\delta_j)})$.
By parts 1 and 8 of Lemma \ref{lgsB}, $\rho=\kappa_\rho$ and $max_1(\rho)=\rho$.
By part 2 of Lemma \ref{lgsD}, this implies $\rho+\kappa_{\rho\cdot\delta_i^-}= \kappa_{\rho+\rho\cdot\delta_i^-}=\kappa_{\rho\cdot\delta_i}$ and $\rho+\kappa_{\rho\cdot\delta_i^-+1}= \kappa_{\rho+\rho\cdot\delta_i^-+1}=\kappa_{\rho\cdot\delta_i+1}$.
By part 2 of Lemma \ref{lisA} again, this implies that $max(X_j)+\rho+\kappa_{\rho\cdot\delta_i^-}=\kappa_{\rho\cdot(\delta_0+\cdots+\delta_i)}$ and $max(X_j)+\rho+\kappa_{\rho\cdot\delta_i^-+1}=\kappa_{\rho\cdot(\delta_0+\cdots+\delta_i)+1}$.
Therefore, 
$$\kappa_{\rho\cdot(\delta_0+\cdots\delta_i)}\leq X_i < \kappa_{\rho\cdot(\delta_0+\cdot+\delta_i)+1}$$
 i.e. $X_i\subseteq I_{\rho\cdot(\delta_0+\cdots+\delta_i)}$.
 
 If $\xi=0$, $X_0\cup\cdots\cup X_n$ satisfies the conclusion of the claim. 
 So, we may assume $0<\xi$. 
 Let $\delta=\delta_0+\cdots+\delta_n$, $X'=X_0\cup\cdots\cup X_n$ and $X=X'\cup\{\kappa_{\rho\cdot\delta+\xi}\}$.
 It suffices to show $X$ is incompressible.
 By part 8 of Lemma \ref{lgsB}, $\kappa_{\rho\cdot\delta+\xi}=max_1(\kappa_{\rho\cdot\delta})+\xi$.
 Since $max(I_{\rho\cdot\delta})=max_1(\kappa_{\rho\cdot\delta})$,  part 9 of Lemma \ref{lisA} implies $max_1(\kappa_{\rho\cdot\delta})$ is in $X'$.
 Let $max_1(\kappa_{\rho\cdot\delta})=\rho\cdot\tau+\epsilon$ where $\epsilon<\rho$. 
 Since $X'$ is closed, $\rho\cdot\tau\in X$.
Since $\kappa_{\rho\cdot\delta+\xi}=\rho\cdot\tau+(\epsilon+\xi)$ and $\epsilon+\xi<\rho$, $X'\cup \{\kappa_{\rho\cdot\delta+\xi}\}$ is closed.
 By part 6 of Lemma \ref{lisA}, $X$ is incompressible.
\qed

\blankline

We remark that the converse of the previous lemma is true.

\begin{lem}
\label{lisG}
If $X$ is a finite nonempty set of ordinals such that $X\subseteq \kappa_\alpha$ for some $\alpha<\theta_1$ then there is an incompressible set $X^+$ such that $X\subseteq X^+$ and $max(index[X^+])=max(index[X])$.
\end{lem}
{\bf Proof.}
There is a finite set $K$ satisfying the hypothesis of the previous lemma such that $index[X]\subseteq K$ and $max(index[X])=max(K)$. 
By the previous lemma, there is an incompressible set $X'$ such that $index[X']=K$. 
By part 3 of Lemma \ref{lisA}, we may assume $\kappa_\alpha\in X'$ whenever $\alpha\in K$.
By part 4 of Lemma \ref{lisA}, $X^+=X\cup X'$ is incompressible.
\qed

\section{The Second Recurrence Theorem for \leqone}
\label{ssecondrone}

Recall that  $\rho$ is an arbitrary additively indecomposable ordinal.

\begin{lem}
\label{lgsE10}
Assume $max_1(\alpha)\not=\infty$ and $[0,\alpha] \not \leq_2 [\alpha+1, max_1(\alpha)]$.
There exists $\lambda<\theta_1$ such that $max_1(\alpha)=\alpha+max_1(\kappa_\lambda)$.
\end{lem}
{\bf Proof.}
Let $\delta$ satisfy $\alpha+\delta=max_1(\alpha)$.
Let $\lambda$ be maximal such that $\kappa_\lambda\leq \delta$. 
We claim $max_1(\kappa_\lambda)=\delta$ from which the conclusion of the lemma follows.
If $ max_1(\kappa_\lambda)<\delta$ then $ \kappa_{\lambda+1}=max_1(\kappa_\lambda)+1\leq \delta$ (by part 8 of Lemma \ref{lgsB}) which contradicts the choice of $\lambda$.
Therefore, $\delta\leq max_1(\kappa_\lambda)$.
Argue by contradiction and assume $\delta<max_1(\kappa_\lambda)$.
This implies $\kappa_\lambda\leqone \delta+1$ and $0<\kappa_\lambda$.
By the Main Structural Lemma, $[1,\delta+1]\cong[\alpha+1,\alpha+\delta+1]$.
This implies $\alpha+\kappa_\lambda\leqone \alpha+\delta+1$.
Since $\alpha+\kappa_\lambda\leq \alpha+\delta=max_1(\alpha)$,  this implies $\alpha\leqone \alpha+\delta+1$ which contradicts the choice of $\delta$.
\qed


\begin{thm} 
\label{tgsFA}
{\bf (Second Recurrence Theorem for \leqone)}
For any ordinal $\beta$,  if $\rho\cdot\omega^\beta<\theta_1$   then 
\begin{itemize}
\item[$(*)$] \hspace{1in} $ max_1(\kappa_{\rho\cdot\omega^\beta})=\kappa_{\rho\cdot\omega^\beta}+max_1(\kappa_\beta)$
\end{itemize}
\end{thm}
{\bf Proof.}
We will prove $(\ast)$ by induction on those $\beta$ with $\rho\cdot\omega^\beta<\theta_1$. 

\halfblankline

{\bf Claim.} Assume $0<\delta_1<\delta_2$, $\delta_2$ is  additively indecomposable and $\rho\cdot\delta_2+1<\theta_1$. If $max_1(\kappa_{\rho\cdot\delta_i})=\kappa_{\rho\cdot\delta_i}+\mu_i$ for $i=1,2$ then $\mu_1<\mu_2$. 

It suffices to show that $\kappa_{\rho\cdot\delta_2}\leqone \kappa_{\rho\cdot\delta_2}+\mu_1+1$. 
For this, suppose $X\subseteq \kappa_{\rho\cdot\delta_2}$ and $Y\subseteq [\kappa_{\rho\cdot\delta_2},\kappa_{\rho\cdot\delta_2}+\mu_1]$ are finite with $X$ and $Y$ closed. 
We will show there exists $\tilde{Y}\subseteq \kappa_{\rho\cdot\delta_2}$ such that $X<\tilde{Y}$ and $X\cup \tilde{Y}$ is a covering of $X\cup Y$.
Notice that $\eta\not\leq_k \zeta$ for $\eta\in X$, $\zeta\in Y$ and $k=1,2$.
By Lemma \ref{lblC}, it suffices to find  $\tilde{Y}\subseteq \kappa_{\rho\cdot\delta_2}$ which is a covering of $Y$ with $X<\tilde{Y}$.

By part 1 of Lemma \ref{lmslA}, $I_{\rho\cdot\delta_1}=[\kappa_{\rho\cdot\delta_1},\kappa_{\rho\cdot\delta_1}+\mu_1]$ is a covering of $[\kappa_{\rho\cdot\delta_2},\kappa_{\rho\cdot\delta_2}+\mu_1]$.
Therefore, $I_{\rho\cdot\delta_1}$ contains a covering of $Y$.
Choose $\gamma<\delta_2$ such that $X<\kappa_{\rho\cdot\gamma}$. 
Since $\delta_2$ is additively indecomposable, $\gamma+\delta_1<\delta_2$.
Since $I_{\rho\cdot(\gamma+\delta_1)}=I_{\rho\cdot\gamma+\rho\cdot\delta_1}\cong I_{\rho\cdot\delta_1}$ by the First Recurrence Theorem for \leqone, $I_{\rho\cdot(\gamma+\delta_1)}$ contains a covering of $Y$.
Moreover, $X<\kappa_{\rho\cdot\gamma}<I_{\rho\cdot(\gamma+\delta_1)}<\kappa_{\rho\cdot\delta_2}$.

\halfblankline

Assume $\rho\cdot\omega^\beta<\theta_1$ and $(\ast)$ holds when $\beta$ is replaced by $\beta'$ for all $\beta'<\beta$. 
Let $\alpha=\rho\cdot\omega^\beta$.
There exists $\lambda<\theta_1$ such that $max_1(\kappa_\alpha)=\kappa_\alpha+max_1(\kappa_\lambda)$.
This is clear if $\alpha+1=\theta_1$ (in which case $max_1(\kappa_\alpha)=\infty$)  and follows from Lemma \ref{lgsE10} otherwise. 
By the induction hypothesis and the claim above, $max_1(\kappa_{\beta'})<max_1(\kappa_\lambda)$ whenever $\beta'<\beta$. 
Therefore, $\beta\leq\lambda$. 

Argue by contradiction and assume $\beta<\lambda$.
This implies $max_1(\kappa_\beta)<max_1(\kappa_\lambda)$.
Therefore, $\kappa_\alpha\leqone \kappa_\alpha+max_1(\kappa_\beta)+1$.

By Lemma \ref{lisG}, there exists  $R$ such that $R$ is incompressible  and $max(Y)=max_1(\kappa_\beta)$.
We may assume $0\in R$.
Let $Y=\kappa_\alpha + R$.
Since $0\in R$, $\kappa_\alpha\in Y$.
Since $\kappa_\alpha\leqone\kappa_\alpha+max_1(\kappa_\beta)+1$, there is a covering $\tilde{Y}$ of $Y$ contained in $ \kappa_\alpha$ with $0<\tilde{Y}$.
Since $\kappa_\alpha\leqone \kappa_\alpha+max_1(\kappa_\beta)$, $\kappa_\alpha\leqone \xi$ for all $\xi\in Y$.
Letting $\mu$ be the least element of $\tilde{Y}$, this implies that $\mu\leqone\xi$ for all $\xi\in\tilde{Y}$ which in turn implies that $\tilde{Y}\subseteq I_{\alpha'}$ for some $\alpha'<\alpha$.
Since $0<\tilde{Y}$, $0<\alpha'$.
Since $\tilde{Y}$ is closed, $\mu$ is divisible by $\rho$. 
By  part 8 of Lemma \ref{lgsB},  $\alpha'$ is divisible by $\rho$.
There are ordinals $\gamma'$ and $\beta'$ such that $\alpha'=\rho\cdot(\gamma'+\omega^{\beta'})$.
 Since $\alpha'<\alpha=\rho\cdot\omega^\beta$, $\beta'<\beta$.
 By the induction hypothesis, $max_1(\kappa_{\alpha'})=\kappa_{\alpha'}+max_1(\kappa_{\beta'})$.
Let $h$ be the covering of $Y$ onto $\tilde{Y}$ and define a function $f$ on $R$ by $h(\kappa_\alpha+\chi)= \kappa_{\alpha'}+f(\xi)$ for $\chi \in Y$.
Since $h$ maps $Y$ into $I_{\alpha'}$, the range of $f$ is contained in $[0,max_1(\kappa_{\beta'})]$.
 By the Main Structural Lemma, $[\kappa_\alpha+1,\infty)\cong [1,\infty)$ and $[\kappa_{\alpha'}+1,\infty)\cong [1,\infty)$. 
 This and the fact that $h$ is a covering easily imply $f$ is a covering.
 SInce $max_1(\kappa_{\alpha'})<\kappa_\alpha$, this contradicts the fact that $R$ is incompressible.
 \qed
 
 \begin{cor}
 For any ordinals $\gamma$ and  $\beta$,  if $\rho\cdot(\gamma+\omega^\beta)<\theta_1$   then 
$$ max_1(\kappa_{\rho\cdot(\gamma+\omega^\beta)})=\kappa_{\rho\cdot(\gamma+\omega^\beta)}+max_1(\kappa_\beta)$$
 \end{cor}
 {\bf Proof.} Combine both the First and Second Recurrence Theorems for \leqone.
 \qed

\begin{cor}
\label{cgsG} 
If $0<\alpha<\theta_1$ then $\kappa_\alpha \leqone \kappa_\alpha + \kappa_\alpha$ iff $\alpha$ is an epsilon number (i.e. $\omega^\alpha = \alpha$) which is greater than $\rho$.
In particular, either $\theta_1=\infty$ or $\theta_1=\theta+1$ where $\theta$ is an epsilon number greater than $\rho$.
\end{cor}
{\bf Proof.}
Assume $\alpha<\theta_1$ is a positive ordinal.

($\Rightarrow$) Assume $\kappa_\alpha \leqone \kappa_\alpha + \kappa_\alpha$. By part 8 of Lemma \ref{lgsB}, $\alpha=\rho\cdot\delta$ for some limit ordinal $\delta$.
Therefore, $\rho<\alpha$.
There are ordinals $\gamma$ and $\beta$ such that $\beta>0$ and $\delta=\gamma+\omega^\beta$.
By the theorem, $\kappa_\alpha+\kappa_\alpha\leq max_1(\kappa_\alpha)=\kappa_\alpha+max_1(\kappa_\beta)$. 
Therefore, $\kappa_\alpha\leq max_1(\kappa_\beta)$ implying $\alpha\leq\beta$. 
This implies 
$$\alpha=\rho\cdot(\gamma+\omega^\beta)\geq \omega^\beta\geq \omega^\alpha$$ 
Therefore, $\alpha=\omega^\alpha$ i.e. $\alpha$ is an epsilon number.

($\Leftarrow$) Assume $\alpha$ is an epsilon number greater than $\rho$. Since $\alpha=\rho\cdot\omega^\alpha$, the theorem implies that $max_1(\kappa_\alpha)= \kappa_\alpha+max_1(\kappa_\alpha)$.
Therefore, $\kappa_\alpha\cdot\omega\leq max_1(\kappa_\alpha)$ implying $\kappa_\alpha\leqone \kappa_\alpha\cdot\omega$.
\qed

\halfblankline

The proof shows that $\kappa_\alpha\leqone \kappa_\alpha\cdot\omega$ when $\alpha<\theta_1$ is an epsilon number greater than $\rho$. 
In fact,  $\kappa_\alpha\leqone\kappa_\alpha\cdot(\omega+1)$ by the Second Recurrence Theorem for $\leqtwo$ from [\ref{Ca??}].
Moreover, $max_1(\kappa_\alpha)=\kappa_\alpha\cdot(\omega+1)$ when $\alpha$ is the least epsilon number greater than $\rho$.

\section{The First Recurrence Theorem for \leqtwo}
\label{sfirstrtwo}

Recall that $\rho$ is an arbitrary additively indecomposable ordinal.

 In the case $\alpha$ is not an epsilon number greater than $\rho$, the first corollary to the Second Recurrence Theorem for \leqone\ provides a description of $I_\alpha$ in terms of the intervals $I_{\alpha'}$ with $\alpha'<\alpha$. We next study the structure of $I_\alpha$  by considering its decomposition into intervals determined by  \leqtwo\ when $\alpha$ is an epsilon number greater than $\rho$.
 
 
\begin{thm}
\label{lgsH}
{\bf (Recurrence Theorem for Small Intervals)}
Assume $\alpha<\theta_1$  is of the form $\rho\cdot\lambda$ where $\lambda$ is infinite and additively indecomposable. Also assume  $\kappa_\alpha$ divides  $\delta\in ORD$.
\begin{enumerate}
 \item $[0,\delta]\not \leq_2 [\delta+1,\delta +\kappa_\alpha)$.
 \item  $ [1,\kappa_\alpha] \cong[\delta +1, \delta+ \kappa_\alpha]$
\end{enumerate}
\end{thm}
{\bf Proof.} We will  prove parts 1 and 2 simultaneously by induction on those $\delta$ which are divisible by $\kappa_\alpha$. 

Assume $\delta$ is divisible by $\kappa_\alpha$ and  conditions 1 and 2 hold for all $\delta'<\delta$ with $\delta'$ divisible by $\kappa_\alpha$. 

For $\delta=0$, both conditions are trivial. 
So, we may assume $\delta>0$. 
By Lemma \ref{lgsDA} and the induction hypothesis, for any finite closed $Y\subseteq \kappa_\alpha$ there are cofinally many subsets of $\delta$ which are isomorphic to $Y$. 

To show conditions 1 and 2 hold, let $J$ be the collection of $\gamma\in[1,\kappa_\alpha]$ such that
\begin{enumerate}
\item[(a)] $[0,\delta]\not\leq_2[\delta+1,\delta+\gamma)$
\item[(b)] $[\delta+1,\delta+\gamma]\cong [1,\gamma]$
\end{enumerate}
We will prove $J=[1,\kappa_\alpha]$ by induction. Suppose $\gamma\in [1, \kappa_\alpha]$ and $\gamma'\in J$ whenever $1\leq\gamma'<\gamma$. 

By the Main Structural Lemma, (b) follows from (a).
To show that (a) holds, argue by contradiction and let $\delta'\leq \delta$ and $\gamma'<\gamma$ with $\delta'\leq_2 \delta+\gamma'$. 
Choose $\alpha'<\alpha$ such that $\gamma'<\kappa_{\alpha'}$. 
By Lemma \ref{lgsF}, there exists incompressible  $Y$ with $max(Y)=\max_1(\kappa_{\alpha'})$. 
In particular, $[0,\gamma)$ does not contain a covering of $Y$.
Since the inductive hypothesis implies that $[1,\gamma)\cong [\delta+1,\delta+\gamma)$, $[\delta+1,\delta+\gamma)$ does not contain a covering of $Y$.
As noted above, there are cofinally many isomorphic copies of $Y$ below $\delta$ and therefore, by part 1 of Lemma \ref{lblA}, there are cofinally $\tilde{Y}$ below $\delta'$ such that $\tilde{Y}$ is a covering of $Y$. 
Since $\delta'\leqtwoup \delta+\gamma$, there are cofinally many coverings of $Y$ below $\delta+\gamma$ -- contradiction.
\qed

\begin{cor}
\label{cgsH1}
Assume $\alpha+1<\theta_1$.
There exist $\delta$ and $\beta$ with  $\beta<\alpha$ such that $max_1(\kappa_\alpha)=\kappa_\alpha\cdot\delta+\max_1(\kappa_\beta)$.
 \end{cor}
{\bf Proof.} 
The corollary follows from part 8 of Lemma \ref{lgsB} if $\alpha$ is not of the form $\rho\cdot\lambda$ where $\lambda$ is an infinite limit ordinal.
So, assume $\alpha=\rho\cdot\lambda$ where $\lambda$ is a limit ordinal. 
The corollary follows from the first corollary to the Second Recurrence Theorem for \leqone\ if $\lambda$ is not additively indecomposable.
So, assume $\lambda$ is additively indecomposable.

Choose $\delta$ divisible by $\kappa_\alpha$ such that  $\kappa_\alpha\leq \delta \leq max_1(\kappa_\alpha)<\delta+\kappa_\alpha$. 
By part 2 of the theorem, $[\delta+1,\delta+\kappa_\alpha)\cong [1,\kappa_\alpha)$.
The corollary now follows by an argument similar to that used for Lemma \ref{lgsE10}.
\qed

\blankline

The first corollary of the Second Recurrence Theorem for \leqone\ provides a computation of the $\beta$ in the corollary above in the case $\alpha$ is not an epsilon number greater than $\rho$. We will not pursue the computation of $\lambda$ further in this paper.


\begin{lem}
\label{lgsI}
Assume $\alpha<\theta_1$ is an epsilon number greater than $\rho$. If $\nu \in I_\alpha$ is divisible by $\kappa_\alpha$ and minimal with respect to $\leqtwo$ then $max_1(\nu)= max_1(\kappa_\alpha)$.
\end{lem}
{\bf Proof.} Assume $\nu\in I_\alpha$ is divisible by $\kappa_\alpha$. 
Since $\nu\in I_\alpha$, $max_1(\nu)\leq max_1(\kappa_\alpha)$.
We will show that   for all $\xi\in I_\alpha$, if $\nu\leq\xi$ then $\nu \leqone \xi$. Argue by induction on $\xi$.

Suppose $\xi \in I_\alpha$, $\nu\leq\xi$ and $\nu\leqone\xi'$  for all $\xi'\in[\nu,\xi)$. 
Since $\nu\leqone\nu$, we may assume $\nu<\xi$.
If $\xi$ is a limit ordinal then $\nu\leqone\xi$ by part 2 of Lemma \ref{lblE}. 
So we may assume $\xi = \eta+1$ where $\eta\geq \nu$.

Suppose $X\subseteq \nu$ and $Y\subseteq [\nu,\xi)=[\nu,\eta]$ are finite closed sets.
Since $\kappa_\alpha \leqone \xi$, there are cofinally many $\tilde{Y}$ below $\kappa_\alpha$ such that $\tilde{Y}$ is a covering of $Y$ by part 1 of Lemma \ref{lblA}. 
By part 2 of the previous lemma, there are cofinally many $\tilde{Y}$ below $\nu$ such that $\tilde{Y}$ is a covering of $Y$. 
Choose such a $\tilde{Y}$ with $X<\tilde{Y}$.
Since $\nu$ is minimal with respect to $\leqtwo$ and $\nu \leqone \eta$, we must have $[0,\nu)\not\leq_2 [\nu,\eta]$.
By Lemma \ref{lblC}, this implies that $X\cup \tilde{Y}$ is a covering of $X\cup Y$.
\qed


\begin{lem}
\label{lgsI1}
Assume $\alpha<\theta_1$ is an epsilon number greater than $\rho$.
The collection of $\nu\in I_\alpha$ such that $\nu$ is divisible by $\kappa_\alpha$ and minimal with respect to \leqtwo\ is topologically closed.
\end{lem}
{\bf Proof.} Let $A$ be a nonempty subset of the collection of ordinals $\nu\in I_\alpha$ which are divisible by $\kappa_\alpha$ and minimal with respect to \leqtwo. 
Let $\mu$ be the least $\xi$ such that $\nu\leq\xi$ for all $\nu\in A$. 
Clearly, $\mu$ is divisible by $\kappa_\alpha$. 
To show $\mu$ is minimal with repect to \leqtwo, argue by contradiction and suppose $\tau<_2\mu$. There exists $\nu\in A$ such that $\tau<\nu\leq\mu$. 
By the previous lemma, $\nu\leqone max_1(\kappa_\alpha)$. 
Therefore, $\nu \leqone \mu$. 
This implies that $\tau <_2 \nu$ -- contradiction.
\qed

\blankline

\begin{dfn}
Assume $\alpha<\theta_1$ is an epsilon number greater than $\rho$.
Define
$$\nu_{\alpha,\xi} \ \ \ (\xi<\theta_2(\alpha)\, ) $$ 
to be the enumeration of the elements of $I_\alpha$ which are divisible by $\kappa_\alpha$ and minimal with respect to \leqtwo.
When $\xi+1<\theta_2(\alpha)$, define $J_{\alpha,\xi}$ to be $[\nu_{\alpha,\xi},\nu_{\alpha,\xi+1})$.
If $\xi+1=\theta_2(\alpha)$, define $J_{\alpha,\xi}$ to be $I_\alpha\cap [\nu_\xi,\infty)$.
Define $\overline{J}_{\alpha,\xi}$ to be the topological closure of $J_{\alpha,\xi}$ i.e. add the least proper upper bound of $J_{\alpha,\xi}$ if it is not already in $J_{\alpha,\xi}$.
\end{dfn}

The enumeration in the definition depends on the parameter $\rho$. Later, when we need to make this dependence clear, we will write $\nu^\rho_{\alpha,\xi}$ for $\nu_{\alpha,\xi}$ and $\theta^\rho_2(\alpha)$ for $\theta_2(\alpha)$. On the other hand, when $\alpha$ is clear from the context, we will write $\nu_\xi$ for $\nu_{\alpha,\xi}$ and $\theta_2$ for $\theta_2(\alpha)$.
We treat the notation $J_{\alpha,\xi}$ and $\overline{J}_{\alpha,\xi}$ similarly.

Notice that when $\xi+1<\theta_2$, $\overline{J}_{\alpha,\xi}=[\nu_\xi,\nu_{\xi+1}]$, and when $\xi+1=\theta_2$, $\overline{J}_{\alpha,\xi}=J_{\alpha,\xi}$.


\begin{lem}
\label{lgsJ}
Assume $\alpha<\theta_1$ is an epsilon number greater than $\rho$. 
\begin{enumerate}
\item
$\nu_0=\kappa_\alpha$ and $max_2(\nu_0)=\nu_0$.
\item
 $\xi\mapsto\nu_\xi$  is continuous.
\item
If $\theta_2\not=\infty$ then there is an ordinal $\theta$ such that $\theta+1=\theta_2$.
\item
If $\xi<\theta_2$ then $[0,\nu_\xi)\not\leq_2[\nu_\xi,\infty)$.
\item 
If $\eta<\xi<\theta_2$ then $J_\eta<J_\xi$ and $J_\eta \not\leq_2 J_\xi$.
\item 
 If $\xi<\theta_2$ then $[\kappa_\alpha,\nu_\xi)=\bigcup_{\eta<\xi}J_\eta$.
 \item
 $I_\alpha=\bigcup_{\eta<\theta_2}J_\eta$.
\item 
If $\xi<\theta_2$ then $max_1(\nu_\xi)= max_1(\kappa_\alpha)$.
 \item 
 If $\xi+1 < \theta_2$ then 
 \begin{center}
  $\nu_\xi \leqtwo \delta$  \ iff \ $\delta \in [\nu_\xi, \nu_{\xi+1})$ and $\delta$ is divisible by $\kappa_\alpha$
  \end{center}
 for all $\delta$.
 \item 
If $\xi+1=\theta_2$ then 
 \begin{center}
  $\nu_{\xi} \leqtwo \delta$  \ iff \ $\delta \in I_\alpha\cap[\nu_\xi, \infty)$ and $\delta$ is divisible by $\kappa_\alpha$
  \end{center}
 for all $\delta$.
\item
     If $\xi+1<\theta_2$ then $\nu_{\xi+1}=max_2(\nu_\xi)+\kappa_\alpha$ and $max_2(\nu_{\xi+1})=\nu_{\xi+1}$.
  \item
    If $\xi<\theta_2$ and $max_2(\nu_\xi)<\infty$ then $[0,max_2(\nu_\xi)]\not\leq_2[max_2(\nu_\xi)+1,\infty)$.
\item
Assume $\xi<\theta_2$,  $\mu$ is divisible by $\kappa_\alpha$, $\chi<\kappa_\alpha$ and $index(\chi)=\gamma$.
 \begin{enumerate}
 \item
 If $\mu+\chi\in J_\xi$ then $\mu+I_\gamma\subseteq J_\xi$.
 \item 
 If $\mu+\chi\in \overline{J}_\xi$ then $\mu+I_\gamma\subseteq \overline{J}_\xi$.
 \end{enumerate}
\item
If $\xi<\theta_2$, $\nu_\xi\leq \mu\in I_\alpha$ and $\mu$ is divisible by $\kappa_\alpha$ then $\nu_\xi\leq_2\uparrow \mu$.
\end{enumerate}
\end{lem}
{\bf Proof.} 
Clearly, $\nu_0=\kappa_\alpha$. 
By part 1 of Lemma \ref{lblD}, $max_2(\kappa_\alpha)=\kappa_\alpha$.
This establishes part 1.

Parts 2 and 3 follow from Lemma \ref{lgsI1}.

Parts 4 through 7 follow immediately from the definitions and Lemma \ref{lgsI1}.

Part 8 is a restatement of Lemma \ref{lgsI}.

\halfblankline

{\bf Claim.} Assume $\delta$ is divisible by $\kappa_\alpha$.
\begin{enumerate}
\item 
If $\delta \leqtwo \delta'$ then $\delta'$ is divisible by $\kappa_\alpha$.
\item 
If $\delta'\leqone \delta$ then $\delta'$ is divisible by $\kappa_\alpha$.
\item
If $\delta\in I_\alpha$ then $\nu_\xi\leqtwo\delta$ for some $\xi<\theta_2$.
\end{enumerate}

Parts 1 and 2 follow from the Recurrence Theorem for Small Intervals.
For part 3, suppose $\delta\in I_\alpha$. 
There exists $\mu\leqtwo\delta$ which is minimal in $\leqtwo$. 
Since $\mu\leqone\delta\in I_\alpha$, $\mu\in I_\alpha$.
By part 2, $\mu$ is divisible by $\kappa_\alpha$.
Therefore, $\mu=\nu_\xi$ for some $\xi<\theta_2$.

\halfblankline

Parts 9 and 10 of the lemma follow easily from  parts 1 and 3 of the claim and part 4 of the lemma.

For part 11, assume $\xi+1<\theta_2$. Part 1 of the claim implies that $max_2(\nu_\xi)$ is divisible by $\kappa_\alpha$. 
Therefore, $max_2(\nu_\xi)+\kappa_\alpha$ is the least ordinal greater than $max_2(\nu_\xi)$ which is divisible by $\kappa_\alpha$.
By part 9 of the lemma, $\nu_{\xi+1}=max_2(\nu_\xi)+\kappa_\alpha$. 

To see $max_2(\nu_{\xi+1})=\nu_{\xi+1}$ argue by contradiction and assume $\nu_{\xi+1}<_2 \delta$.
By part 1 of Lemma \ref{lblD}, there exists $\gamma$ such that $\gamma<_1\nu_{\xi+1}$ and $max_2(\nu_\xi)<\gamma$.
By part 2 of the claim, $\gamma$ is divisible by $\kappa_\alpha$ which contradicts $\nu_{\xi+1}=max_2(\nu_\xi)+\kappa_\alpha$.

For part 12, assume $\xi<\theta_2$ and $max_2(\nu_\xi)\not=\infty$.
Argue by contradiction and assume $\nu\leq max_2(\nu_\xi)<\delta$ and $\nu\leqtwo\delta$.
This implies $\nu\leqone max_2(\nu_\xi)$. 
Clearly, $\nu,\delta\in I_\alpha$.
Part 1 of the claim implies $max_2(\nu_\xi)$ is divisible by $\kappa_\alpha$. 
By part 2 of the claim, $\nu$ is divisible by $\kappa_\alpha$.
By part 1 of the claim, $\delta$ is divisible by $\kappa_\alpha$.
By part 10 of the lemma, $\xi+1<\theta_2$.
By part 11, $\nu_{\xi+1}\leq\delta$ contradicting part 4.

For part 13, assume $\xi<\theta_2$, $\mu$ is divisible by $\kappa_\alpha$, $\chi<\kappa_\alpha$ and $\gamma$ is the index of $\chi$. 

To establish part (a), assume $\mu+\chi\in J_\xi$.
Since the least element of $J_\xi$ is $\nu_\xi$, $\nu_\xi\leq \mu\leq \mu+I_\gamma$. 

First suppose $\xi+1<\theta_2$. 
In this case, $J_\xi=[\nu_\xi,\nu_{\xi+1})$ and $\mu<\nu_{\xi+1}$. 
Since $\nu_{\xi+1}$ is divisible by $\kappa_\alpha$, $\mu+I_\gamma<\mu+\kappa_{\gamma+1}<\mu+\kappa_\alpha\leq \nu_{\xi+1}$.

Now suppose $\xi+1=\theta_2$. 
If $\alpha+1=\theta_1$ then $J_\xi=[\nu_\xi,\infty)$ implying $\mu+I_\gamma\subseteq J_\xi$. 
So, we may assume $\alpha+1<\theta_1$. 
By Corollary \ref{cgsH1}, $I_\alpha=[\kappa_\alpha, \kappa_\alpha\cdot\delta+max_1(\kappa_\beta))$ for some $\delta$ and $\beta<\alpha$.
Therefore, $J_\xi=[\nu_\xi, \kappa_\alpha\cdot\delta +max_1(\kappa_\beta)]$.
If $\mu+\chi< \kappa_\alpha\cdot\delta$ the argument for the case $\xi+1<\theta_2$ shows $\mu+I_\gamma\subseteq J_\xi$.
So, we may assume $\kappa_\alpha\cdot\delta\leq\mu+\chi$.
In this case $\kappa_\alpha\cdot\delta=\mu$ and $\gamma\leq\beta$ easily implying 
$\mu+I_\gamma\subseteq J_\xi$.

To establish part (b), assume $\mu+\chi\in \overline{J}_\xi$. 
If $\mu+\chi\in J_\xi$ then $\mu+I_\gamma\subseteq J_\xi\subseteq \overline{J}_\xi$ by part (a).
So, we may assume that $\mu+\chi\not\in J_\xi$.
This implies $\xi+1<\theta_2$ and $\mu+\chi=\nu_{\xi+1}$.
Since $\nu_{\xi+1}$ is divisible by $\kappa_\alpha$, $\mu=\nu_{\xi+1}$, $\chi=0$ and $\gamma=0$.
Therefore, $\mu+I_\gamma=\mu+I_0=\mu+\{0\}=\{\mu\}\subseteq \overline{J}_\xi$.

For part 14, assume $\mu$ is a multiple of $\kappa_\alpha$ with $\nu_\xi\leq\mu\in I_\alpha$.
To show $\nu_\xi\leq^\rho_2\uparrow \mu$, assume $X$ is a finite $\rho$-closed subset of $\nu_\xi$ and $\cal Y$ is a collection of finite $\rho$-closed subsets of $\nu_\xi$ which is cofinal in $\nu_\xi$ such that $X<Y$ whenever $Y\in \cal Y$ and $X\cup Y_1\cong_\rho X\cup Y_2$ whenever $Y_1,Y_2\in \cal Y$. 

Assume $Y\in \cal Y$.
We will show there exist cofinally many finite closed subset $\tilde{Y}$ of $\mu$ such that  $X\cup \tilde{Y}$ is a covering of $X\cup Y$. 
By part 8, $\nu_\xi\leqone \mu$.
By Lemma \ref{lblC1}, we may assume that for each $\sigma\in X$  there exists $\tau\in Y$ such that $\sigma \leq_2\tau$. 
By part 3 of Lemma \ref{lblE}, this assumption implies that if $\sigma\in X$ then $\sigma\leq_2 \nu_\xi$.
Since $\nu_\xi$ is minimal in \leqtwo, $X=\emptyset$. 
Fix $Y\in \cal Y$. 
Since $\kappa_\alpha\leq_1 \nu_\xi$, there are cofinally many covering of $Y$ below $\kappa_\alpha$.
Since $\mu$ is a multiple of $\kappa_\alpha$, RTSI implies there are cofinally many coverings of $Y$ below $\mu$.
\qed


\begin{lem}
\label{lrec2B}
Assume $\alpha<\theta_1$ is an epsilon number greater than $\rho$.
The following hold for  $\xi+1<\theta_2$ and $\eta<\theta_2$.
\begin{enumerate}
\item
$I_\alpha  \cong I_\alpha \cap [\nu_{\xi+1},\infty) $.
\item
$\xi+1+\eta<\theta_2$
\item
$\nu_{\xi+1+\eta}=max_2(\nu_{\xi})+\nu_\eta$
\item
$max_2(\nu_{\xi+1+\eta})=max_2(\nu_\xi)+max_2(\nu_\eta)$
\end{enumerate}
\end{lem}
{\bf Proof.}
By part 11 of the previous lemma and the Main Structural Lemma, the map $h(\gamma)=max_2(\nu_\xi)+\gamma$ is an isomorphism of $[1,\infty)$ and $[max_2(\nu_\xi)+1,\infty)$.

\halfblankline

{\bf Claim.} 
\begin{enumerate}
\item
For $\beta>0$, $\beta$ is divisible by $\kappa_\alpha$ iff $h(\beta)$ is divisible by $\kappa_\alpha$.
\item
For $\beta>0$ and $k=1,2$, 
\begin{center}
$max_k(\beta)=\infty$ iff $\max_k(h(\beta))=\infty$
\end{center}
and 
$$h(max_k(\beta))=max_k(h(\beta))$$
if $max_k(\beta)\not=\infty$.
\item
$h(\kappa_\alpha)=\nu_{\xi+1}$
\item
If $max_1(\kappa_\alpha)\not=\infty$ then $h(max_1(\kappa_\alpha))=max_1(\kappa_\alpha)$.
\item
For $\beta>0$, $\beta\in I_\alpha$ iff $h(\beta)\in I_\alpha$.
\item
For $\beta>0$, $\beta$ is minimimal with respect to \leqtwo\ iff $h(\beta)$ is minimal with respect to \leqtwo.
\item
For $\zeta<\theta_2$, $h(\nu_\zeta)=\nu_{\xi+1+\zeta}$
\end{enumerate}

By part 9 of the previous lemma, $max_2(\nu_\xi)$ is divisible by $\kappa_\alpha$. Part 1 of the claim follows (alternatively, one could use Lemma \ref{lgsD0} which says $\kappa_\alpha$ is additively indecomposable). 

Part 2 follows from the fact $h$ is an isomorphism.

For part 3, notice
 $h(\kappa_\alpha)=max_2(\nu_\xi)+\kappa_\alpha=\nu_{\xi+1}$
by part 11 of the previous lemma.

Since $max_1(\nu_{\xi+1})=max_1(\kappa_\alpha)$ by part 8 of the previous lemma,
part 4 follows from
$$h(max_1(\kappa_\alpha))=max_1(h(\kappa_\alpha))=max_1(max_2(\nu_\xi)+\kappa_\alpha)=max_1(\nu_{\xi+1})$$

Part 5 follows from parts 3 and 4.

For part 6, suppose $\beta>0$.
First notice $\beta$ is minimal with respect to \leqtwo\  iff $\beta$ is minimal with respect to \leqtwo\ in $[1,\infty)$.
Since $h$ is an isomorphism, $\beta$ is minimal with respect to \leqtwo\ in $[1,\infty)$ iff $h(\beta)$ is minimal with respect to \leqtwo\ in $[max_2(\nu_\xi)+1,\infty)$.
By part 12 of the previous lemma,  $h(\beta)$ is minimal with respect to \leqtwo\ in  $[max_2(\nu_\xi)+1,\infty)$ iff $h(\beta)$ is minimal with respect to \leqtwo.

By parts 1 and 3 through 6, $h(\nu_\zeta)$ $(\zeta<\theta_2)$ enumerates the ordinals which are divisible by $\kappa_\alpha$, minimal with respect to \leqtwo, in $I_\alpha$ and at least $\nu_{\xi+1}$. Since $\nu_{\xi+1+\zeta}$ $(\xi+1+\zeta<\theta_2)$ enumerates the same family of ordinals, $h(\nu_\zeta)=\nu_{\xi+1+\zeta}$ for $\zeta<\theta_2$.

This completes the proof of the claim.

\halfblankline

Part 1 of the lemma follows from parts 3 and 4 of the claim and the fact $h$ is an isomorphism.

Parts 2 and 3 of the lemma follow from part 7 of the claim.

The following establishes part 4 of the lemma.
\begin{tabbing}
\hspace{.5in} $max_2(\nu_{\xi+1+\eta})$ \ \= = \ $max_2(h(\nu_\eta))$ \hspace{1in} \= (part 7 of the claim) \\
\> = \ $h(max_2(\nu_\eta))$ \> (part 2 of the claim) \\
\>  = \ $max_2(\nu_\xi)+max_2(\nu_\eta)$ \> (definition of $h$)
\end{tabbing}
\qed

\begin{thm}
\label{trec2C}
Assume $\alpha<\theta_1$ is an epsilon number greater than $\rho$.
Either $\theta_2=\infty$ or $\theta_2=\theta+1$ for some infinite additively indecomposable ordinal $\theta$.
\end{thm}
{\bf Proof.}
Assume $\theta_2\not=\infty$. 
By part 3 of Lemma \ref{lgsJ}, there exists an ordinal $\theta$ such that $\theta_2=\theta+1$.
By Corollary \ref{cgsG}, $\kappa_\alpha\leqone\kappa_\alpha+\kappa_\alpha$.
By part 1 of Lemma \ref{lblD}, $\kappa_\alpha\not\leq_1\kappa_\alpha+\kappa_\alpha$.
By part 7 of Lemma \ref{lgsJ}, $\theta\not=0$   (and $\nu_1=\kappa_\alpha+\kappa_\alpha$ by parts 1 and 8 of Lemma \ref{lgsJ}).
By part 2 of the previous lemma, $\xi+1+\theta=\theta$ whenever $\xi+1\leq \theta$.
Since $\theta\not= 0$, this implies that $\theta$ is infinite and additively indecomposable.
\qed


\begin{thm}
{\bf (First Recurrence Theorem for \leqtwo)}
Assume $\alpha<\theta_1$ is an epsilon number greater than $\rho$.
If $\xi+1<\theta_2$  and $0<\eta<\theta_2$ then $\overline{J}_{\xi+\eta}\cong \overline{J}_\eta$.  Moreover, $\overline{J}_0\cong \overline{J}_1$.
\end{thm}
{\bf Proof.}
As in the proof of the Lemma \ref{lrec2B}, part 11 of Lemma \ref{lgsJ} and the Main Structural Lemma imply the map $h(\gamma)=max_2(\nu_\xi)+\gamma$ is an isomorphism of $[1,\infty)$ and $[max_2(\nu_\xi)+1,\infty)$. 
 If $\eta+1=\theta_2$ then $\xi+\eta=\eta$ by the previous theorem from which the conclusion follows trivially.
So, we may assume that $\eta+1<\theta_2$.
By parts 3 and 4 of Lemma \ref{lrec2B}, $\overline{J}_\zeta\cong \overline{J}_{\xi+1+\zeta}$ for all $\zeta<\theta_2$. 
When $\xi=0$, we have $\overline{J}_\zeta\cong \overline{J}_{1+\zeta}$ for all $\zeta<\theta_2$.
Choosing $\zeta=0$,   $\overline{J}_0\cong \overline{J}_1$. 
Choosing $\zeta$ such that  $\eta=1+\zeta$, 
$$\overline{J}_{\xi+\eta}=\overline{J}_{\xi+1+\zeta}\cong \overline{J}_{\zeta}\cong \overline{J}_{1+\zeta}=\overline{J}_\eta$$
\qed



\begin{center}
 REFERENCES
\end{center}

\begin{enumerate}
\item
\label{Ba75}
 J. Barwise, {\bf Admissible Sets and Structures}, 
Springer-Verlag, Berlin, 1975.
\item
\label{Be85}
 M.J. Beeson, {\bf Foundations of Constructive Mathematics},
Spinger-Verlag, Berlin, 1985.
\item
\label{Ca09}
T. Carlson, {\it Patterns of resemblance of order 2}, Annals of Pure and Applied Logic 158 (2009), pp. 90-124.
\item
\label{Ca16}
T. Carlson, {\it Generalizing Kruskal's Theorem to pairs of cohabitating trees}, Archive for Mathematical Logic 55 (2016), pp. 37-48.
\item
\label{Ca??}
T. Carlson, {\it Structural Properties of \calRtwo\ Part II}, preprint.
\end{enumerate}

\end{document}